\documentclass{amsart}

\usepackage{amsmath}
\usepackage{amsfonts}
\usepackage{amssymb,enumerate}
\usepackage{amsthm}
\usepackage[all]{xy}

\newtheorem{lem}{Lemma}[section]
\newtheorem{cor}[lem]{Corollary}
\newtheorem{prop}[lem]{Proposition}
\newtheorem{thm}[lem]{Theorem}

\newtheorem{Defn}[lem]{Definition}
\newtheorem{Ex}[lem]{Example}
\newtheorem{Question}[lem]{Question}
\newtheorem{Property}[lem]{Property}
\newtheorem{Properties}[lem]{Properties}
\newtheorem{Discussion}[lem]{Remark}
\newtheorem{Construction}[lem]{Construction}
\newtheorem{Subprops}{}[lem]
\newtheorem{Para}[lem]{}

\newenvironment{ex}{\begin{Ex}\rm}{\end{Ex}}
\newenvironment{question}{\begin{Question}\rm}{\end{Question}}

\newenvironment{para}{\begin{Para}\rm}{\end{Para}}
\newenvironment{disc}{\begin{Discussion}\rm}{\end{Discussion}}

\newcommand{\y}{\mathbf{y}}
\newcommand{\cbc}[2]{#1(#2)}

\newcommand{\comp}[1]{\widehat{#1}}
\newcommand{\ideal}[1]{\mathfrak{#1}}
\newcommand{\m}{\ideal{m}}
\newcommand{\n}{\ideal{n}}
\newcommand{\p}{\ideal{p}}
\newcommand{\q}{\ideal{q}}

\newcommand{\pdim}{\operatorname{pd}}	
\newcommand{\supp}{\operatorname{Supp}}	
\newcommand{\pd}{\operatorname{pd}}	
\newcommand{\gdim}{\mathrm{G\text{-}dim}}	
\newcommand{\gkdim}[1]{\mathrm{G}_{#1}\text{-}\dim}	
\newcommand{\ext}{\operatorname{Ext}}

\newcommand{\depth}{\operatorname{depth}}	
\newcommand{\rank}{\operatorname{rank}}	
\newcommand{\Ker}{\operatorname{Ker}}	
\newcommand{\zz}{\mathbb{Z}}

\newcommand{\injdim}{\operatorname{id}}	
\newcommand{\rhom}{\mathbf{R}\mathrm{Hom}}	
\newcommand{\lotimes}{\otimes^{\mathbf{L}}}
\newcommand{\amp}{\operatorname{amp}}
\newcommand{\HH}{\operatorname{H}}
\newcommand{\Hom}{\operatorname{Hom}}	

\newcommand{\coker}{\operatorname{Coker}}
\newcommand{\fd}{\operatorname{fd}}

\newcommand{\spec}{\operatorname{Spec}}

\newcommand{\s}{\mathfrak{S}}

\newcommand{\vf}{\varphi}

\newcommand{\D}{\mathsf{D}}

\newcommand{\len}{\operatorname{length}}

\newcommand{\xra}{\xrightarrow}

\newcommand{\shift}{\mathsf{\Sigma}}

\newcommand{\Pic}{\operatorname{Pic}}

\newcommand{\mspec}{\mathrm{m\text{-}Spec}}
\newcommand{\grade}{\operatorname{grade}}
\newcommand{\cone}{\operatorname{cone}}
\newcommand{\image}{\operatorname{Im}}

\begin{document}

\bibliographystyle{amsplain}

\author{Anders Frankild}
\thanks{A.F., University of Copenhagen, Institute for Mathematical 
Sciences, Department of Mathematics, 
Universitetsparken 5, 2100 K\o benhavn, Denmark, \texttt{frankild@math.ku.dk}}

\author{Sean Sather-Wagstaff}
\thanks{S.S.-W., Department of Mathematics, California State University, Dominguez Hills, 
1000 E.~Victoria St., Carson, CA 90747 USA, \texttt{ssather@csudh.edu}}
\thanks{This research 
was conducted while 
A.F.~was funded by the Lundbeck Foundation and by
Augustinus Fonden, and S.S.-W.~was an NSF Mathematical Sciences 
Postdoctoral Research Fellow and a visitor at the University of 
Nebraska-Lincoln.}

\title[Reflexivity and ring homomorphisms]{Reflexivity and ring homomorphisms 
of finite flat dimension}
\keywords{semidualizing complexes, semidualizing modules,
Gorenstein dimensions, 
G-dimensions, ring homomorphisms} 
\subjclass[2000]{13C13, 13D05, 13D25, 13H10}

\dedicatory{Dedicated to the memory of Saunders Mac Lane.}

\begin{abstract}
In this paper we present a systematic study of the reflexivity 
properties of homologically finite complexes with respect to 
semidualizing complexes 
in the setting of nonlocal rings.  One primary 
focus is the descent of these properties over ring homomorphisms 
of finite flat dimension, presented in terms of 
inequalities between generalized G-dimensions.
Most of these results are new even when the 
ring homomorphism is local.  The main tool for these analyses
is a nonlocal version of the amplitude inequality of Iversen, 
Foxby, and Iyengar.  We provide numerous examples demonstrating the 
need for certain hypotheses and the strictness of many inequalities.
\end{abstract}
\maketitle

\section*{Introduction}

Grothendieck and 
Hartshorne~\cite{hartshorne:rad,hartshorne:lc} introduced the notion 
of a dualizing complex as a tool for understanding cohomology 
theories in algebraic geometry and commutative algebra.
The homological properties of 
these objects and the good behavior of rings admitting them are 
well-documented and of continuing interest and application in these 
fields.

Semidualizing complexes arise in several contexts in 
commutative algebra as natural generalizations of dualizing 
complexes; see~\ref{homothety}.
A dualizing complex for $R$ is semidualizing, as is a
free $R$-module of rank 1.  
Such objects were introduced and studied in the abstract
by Foxby~\cite{foxby:gmarm} and Golod~\cite{golod:gdagpi}
in the case where $C$ is a module.  
The investigation of the general situation begins with the work of 
Christensen~\cite{christensen:scatac} and continues with, e.g., 
\cite{takahashi:hiatsb,frankild:sdcms,gerko:sdc,gerko:osmagi,sather:divisor}.   

The utility of these complexes was first demonstrated in the work of Avramov 
and Foxby~\cite{avramov:rhafgd} where the dualizing complex $D^{\vf}$
of a local ring
homomorphism $\vf\colon R\to S$ of finite flat dimension (or more 
generally of finite 
G-dimension) is used as one way to relate the Bass series of $R$ to 
that of $S$;  see~\ref{paraBass}.  When $\vf$ is module-finite, its dualizing 
complex is $\rhom_R(S,R)$, which is semidualizing for $S$.  (For the 
general case, see~\cite{avramov:rhafgd}.)
This provides another 
generalization of dualizing complexes: if $R$ is 
Gorenstein, then $D^{\vf}$ is dualizing for $S$.
It is believed 
that $D^{\vf}$ will give insight into the 
so-called 
composition 
question for homomorphisms of finite G-dimension.

A semidualizing complex $C$ gives rise to the 
category of $C$-reflexive complexes 
equivalently, the category of complexes
of finite $\text{G}_C$-dimension;  see~\ref{gd}.
When $C$ is dualizing, every homologically finite complex 
$X$ is $C$-reflexive~\cite{hartshorne:rad}.  
On the other hand, a 
complex is $R$-reflexive exactly when it has finite 
G-dimension
as defined by Auslander and Bridger~\cite{auslander:adgeteac,auslander:smt} for modules
and Yassemi~\cite{yassemi:gd} for complexes.  
This notion was introduced and studied in general
by Foxby~\cite{foxby:gmarm} and Golod~\cite{golod:gdagpi}
when $C$ and $X$ are modules, and
by Christensen~\cite{christensen:scatac} in this generality.

The current paper is part of our ongoing effort 
to increase the understanding of the 
semidualizing complexes and their corresponding reflexive 
complexes. 
More of our work in this direction is found 
in~\cite{frankild:sdcms,sather:divisor} where we forward 
two new perspectives for this study.  
In~\cite{frankild:sdcms} we endow the set 
of shift-isomorphism classes of semidualizing $R$-complexes
with a nontrivial 
metric.  
In~\cite{sather:divisor} S.S.-W.~investigates the 
consequences of the observation that, when $R$ is a normal domain, 
the set 
of isomorphism classes of semidualizing $R$-modules
is naturally a subset of the divisor class group 
of $R$.  
Each of these works 
relies 
heavily on the homological tools developed in the current paper, 
which fall into roughly three categories.

First, we extend a number of results 
in~\cite{christensen:scatac}
from the setting of local rings and local ring 
homomorphisms to the nonlocal realm.  This process is begun in 
Section~\ref{sec5} with an investigation of the behavior of these objects 
under localization, and it is continued in Section~\ref{sec3} where global 
statements are proved over a single ring.

The second advancement in this paper is found in the descent results 
which populate Sections~\ref{sec2}--\ref{sec6}.  
Based in part on the ideas of 
Iyengar and S.S.-W.~\cite{iyengar:golh}, we exploit the 
amplitude inequality of Iversen~\cite{iversen:aifc} and Foxby and 
Iyengar~\cite{foxby:daafuc} in order 
to prove converses of a number of results from~\cite{christensen:scatac}.
These results deal with
the interactions between, on the one hand, semidualizing and reflexive 
complexes, and on the other hand, complexes and ring homomorphisms of 
finite flat dimension.
Most of the results from~\cite{christensen:scatac}
that we focus on are stated there in 
the local setting, and the converses are new even there.  However, our 
work in the earlier sections along with a nonlocal version of the 
amplitude inequality 
extend these converses and 
the original results to the global arena.
Our version of the amplitude inequality is Theorem~\ref{prop201}, 
wherein $\inf(X)$ and $\sup(X)$ are the infimum and supremum, 
respectively of the set $\{i\in\mathbb{Z}\mid\HH_i(X)\neq 0\}$ and 
$\amp(X)=\sup(X)-\inf(X)$.

\medskip

\noindent \textbf{Theorem I.}  
\emph{Let $\vf\colon R\to S$ be a ring homomorphism and
$P$ a homologically 
finite $S$-complex 
with $\fd_R(P)$ finite and
such that $\vf^*(\supp_S(P))$ 
contains $\mspec(R)$.
For each homologically degreewise finite $R$-complex
$X$
there are inequalities
\begin{align*}
\inf(X\lotimes_R P)&\leq\inf(X)+\sup(P)\\
\sup(X\lotimes_R P)&\geq\sup(X)+\inf(P)\\
\amp(X\lotimes_R P)&\geq\amp(X)-\amp(P).
\end{align*}
In particular, 
\begin{enumerate}[\quad\rm(a)]
\item 
$X\simeq 0$ if and only if $X\lotimes_R P\simeq 0$; 
\item 
$X$ is homologically bounded if and only if $X\lotimes_R P$ is so;
\item 
If $\amp(P)=0$, e.g., if $P=S$, then $\inf(X\lotimes_R 
P)=\inf(X)+\inf(P)$.
\end{enumerate}}

\medskip

Section~\ref{sec2} deals for the most part with the behavior of the 
semidualizing and reflexive properties 
with respect to the derived functor $-\lotimes_R 
S$ where $\vf\colon R\to S$ is a ring homomorphism of finite flat 
dimension, that is, with $\fd_R(S)<\infty$.  
As a sample, here is a summary of Theorems~\ref{lem01a1},
\ref{lem01a3}, and~\ref{lem01a4}.

\medskip

\noindent \textbf{Theorem II.}  
\emph{Let $\vf\colon R\to S$ be a ring homomorphism of finite flat 
dimension and $C,C',X$ homologically degreewise finite $R$-complexes.
Assume that every maximal ideal of $R$ is contracted from $S$.
\begin{enumerate}[\quad\rm(a)]
\item \label{itemIa}
The complex $C\lotimes_R S$ is $S$-semidualizing  if and only if 
$C$ is $R$-semidualizing.
\item \label{itemIb}
When $C$ is semidualizing for $R$, there is an equality
$$\gkdim{C}_R(X)=\gkdim{C\lotimes_R S}_S(X\lotimes_R S).
$$ 
In particular, 
$X\lotimes_R S$
is $C\lotimes_R S$-reflexive if and only if $X$ is $C$-reflexive.
\item\label{itemIc}
If the induced map on Picard groups $\Pic(R)\to\Pic(S)$ is injective
and $C,C'$ are semidualizing $R$-complexes, then
$C\lotimes_R S$ is isomorphic to $C'\lotimes_R S$ in $\D(S)$
if and only if $C$ is isomorphic to $C'$ in $\D(R)$.
\end{enumerate}
}

\medskip

Section~\ref{sec7} is similarly devoted to the functor $\rhom_R(S,-)$
when $\vf\colon R\to S$ is module-finite.  The 
version of Theorem II for this context is contained in Theorems~\ref{lem102},
\ref{lem105}, and~\ref{lem206}.
We highlight here the characterization of
reflexivity of $\rhom_R(S,X)$ with respect to 
$C\lotimes_R S$ which is in Theorem~\ref{lem107}.  

\medskip

\noindent \textbf{Theorem III.}  
\emph{Let $C,X$ be homologically finite $R$-complexes
with $C$ semidualizing.
If $\vf$ is module-finite with $\fd(\vf)<\infty$
and $\mspec(R)\subseteq\image(\vf^*)$, then
\begin{align*}
\gkdim{C}_R(X)-\pd_R(S)
&\leq\gkdim{\rhom_R(S,C)}_S(X\lotimes_R S)\\
&\leq\gkdim{C}_R(X)+\pd_R(S).
\end{align*}
Thus, 
$X\lotimes_R S$
is $\rhom_R(S,C)$-reflexive if and only if $X$ is $C$-reflexive.
If 
$R$ is local or
$\amp(C)=0=\amp(\rhom_R(S,R))$, then
$$\gkdim{\rhom_R(S,C)}_S(X\lotimes_R S)=\gkdim{C}_R(X).$$
}

In Section~\ref{sec6} we extend results of Section~\ref{sec7} to the
case where $\vf$ is local and admits a Gorenstein factorization
$R\to R' \to S$;  see~\ref{factor}.  To this end, we
use a shift of  the functor 
$\rhom_{R'}(S,-\lotimes_R R')$ in place of
$\rhom_R(S,-)$.  We prove in Theorem~\ref{thm5} that this is 
independent of the choice of Gorenstein factorization and, when $\vf$ is 
module-finite, agrees with $\rhom_R(S,-)$.  The remainder of the 
section is spent documenting the translations of the results from 
Section~\ref{sec7} to this context.

The third focus of this paper is found in the
numerous examples within the text
demonstrating that our results are, in a sense, optimal.
These examples may be of independant interest, as the number of 
explicit computations in this area is somewhat limited.  For this reason, and for
ease of reference, we provide a resume of the more delicate examples 
here.  Note that some of the rings constructed have connected prime spectra, 
and this makes the constructions a tad technical.  
We have taken this approach because 
rings with connected spectra can exhibit particularly nice local-global behavior
and we wanted to make the point that the exemplified behavior 
can occur even when the spectra are connected.

Example~\ref{ex02} 
shows that 
one can have inequalities
$\gkdim{C}_R(X)<\sup(X)$
and
$\gkdim{C}_R(X)<\gkdim{C_{\p}}_{R_{\p}}(X_{\p})$, even when $R$ is local.
Thus, 
$\gkdim{C}_R(X)$ 
cannot be computed as 
the length of a resolution of $X$, 
and the assumption $\amp(C)=0$ is necessary in Lemma~\ref{lem08} and
in the final statement of Lemma~\ref{lem05aa}.

Example~\ref{ex302} provides a
surjective ring homomorphism 
of finite flat dimension that is 
Cohen-Macaulay with nonconstant grade.
Thus, $\spec(S)$ must be connected
in Corollary~\ref{prop203}

Example~\ref{ex102}
shows that one can
have 
$\amp(C)>0$ when $\amp(C_{\m})=0$ for each 
maximal ideal $\m$.
Furthermore, if $C'$ is $C$-reflexive, the inequality
$\amp(C)\leq\amp(C')$ from Corollary~\ref{cor301}
can be strict, even 
when $\amp(C)=0$.  
Thus, the connectedness of $\spec(R)$ is needed in 
Proposition~\ref{prop101} and in
Corollaries~\ref{prop101a} and~\ref{cor301}.

Example~\ref{ex103} provides a ring $R$ with $\spec(R)$ connected where 
\begin{align*}
\inf (C)-\sup (C')&=\inf(\rhom_R(C',C)) < \inf (C)-\inf (C')\\
\inf (C')&< \gkdim{C}_R(C')=\sup (C')\\
\gkdim{A}_R(B) -\sup(B) 
&=
\gkdim{B^{\dagger_C}}_R(A^{\dagger_C})-\inf(C)+\inf(A)\\
& <\gkdim{A}_R(B) -\inf(B)  
\end{align*}
showing that inequalities in Lemma~\ref{lem101} and Proposition~\ref{propDD}
can be strict or not. 

Example~\ref{ex106} shows that strictness can occur in each of the inequalities
\begin{gather*}
\gkdim{C}_R(S)\leq\sup\{\gkdim{C_{\m}}_{R_{\m}}(S_{\m})\mid\m\in\mspec(R)\} \\
\gkdim{C}_R(S)\leq\pd_R(S) \qquad\qquad
\inf(C)\leq\inf(C\lotimes_R S) \\
\inf(C)-\pd_R(S)\leq\inf(\rhom_R(S,C)) \\
\gkdim{C}_R(S)\leq\gkdim{\rhom_R(S,C)}_S(S)+\pd_R(S)
\end{gather*}
from Propositions~\ref{lem06b},
\ref{lem104},
and~\ref{lem201} and from 
Theorems~\ref{lem01a1}
and~\ref{lem102}. 

Example~\ref{ex301} pertains to
Theorems~\ref{lem204}, 
\ref{lem01a3}, 
\ref{lem205}, 
\ref{lem105}, 
and~\ref{lem107}, showing that
one can have 
$\gkdim{C}_R(X)=\infty$ even though 
each of the following is finite:
$\gkdim{C}_R(\rhom_R(S,X))$, 
$\gkdim{\rhom_R(S,C)}_R(\rhom_R(S,X))$,
$\gkdim{C}_R(X\lotimes_R S)$, 
$\gkdim{C\lotimes_R S}_R(X\lotimes_R S)$,
$\gkdim{\rhom_R(S,C)}_R(X\lotimes_R S)$.
Hence, the hypothesis on $\mspec(R)$ is necessary for each result.

As this introduction suggests, most of the results of this paper are 
stated and proved in the framework of the derived category.  We 
collect basic definitions and notations for the reader's convenience 
in Section~\ref{sec1}.

\section{Complexes and ring homomorphisms} \label{sec1}

\noindent\emph{Throughout this work, $R$ and $S$ are commutative Noetherian 
rings and $\vf\colon R\to S$ is a ring homomorphism.}

This section consists of background material and includes most of the 
definitions and notational conventions used throughout the rest of 
this work.

\begin{para} We work 
in the derived category $\D(R)$ whose objects are the $R$-complexes, 
indexed homologically;
references on the subject
include~\cite{gelfand:moha,hartshorne:rad,neeman:tc,verdier:cd,verdier:1}.
For $R$-complexes $X$ and $Y$ 
the left derived tensor product complex 
is denoted $X\lotimes_R Y$ and 
the right derived homomorphism complex is 
$\rhom_R(X,Y)$.  For an integer $n$, 
the $n$th \emph{shift} or \emph{suspension} of $X$ is denoted
$\shift^n X$ where $(\shift^n 
X)_i=X_{i-n}$ and $\partial_i^{\shift^n X}=(-1)^n\partial_{i-n}^X$.
The symbol ``$\simeq$'' indicates an
isomorphism in $\D(R)$ and ``$\sim$'' indicates an isomorphism up to shift.

The \emph{infimum} and \emph{supremum} 
of a complex 
$X$, denoted $\inf(X)$ and $\sup(X)$,
are the infimum and supremum, respectively, of the set 
$\{i\in\zz\mid\HH_i(X)\neq 0\}$, and
the \emph{amplitude} of $X$ is the difference 
$\amp(X)=\sup(X)-\inf(X)$.
The complex $X$ 
is \emph{homologically finite}, respectively 
\emph{homologically degreewise finite}, if its total 
homology module $\HH(X)$, respectively each individual homology 
module $\HH_i(X)$, is a finite $R$-module.
It is \emph{homologically bounded above}, respectively
\emph{homologically bounded below} or 
\emph{homologically bounded}, if $\sup(X)<\infty$, respectively
$\inf(X)>-\infty$ or $\amp(X)<\infty$.
The projective, injective, and flat dimensions of $X$ are denoted 
$\pd_R(X)$, $\injdim_R(X)$, and $\fd_R(X)$, respectively;  
see Avramov and Foxby~\cite{avramov:hdouc}.  
\end{para}

The main objects of study in this paper are the
semidualizing complexes and their reflexive objects, introduced by Foxby~\cite{foxby:gmarm},
Golod~\cite{golod:gdagpi}, and 
Christensen~\cite{christensen:scatac}.  

\begin{para} \label{homothety}
A homologically finite $R$-complex $C$ 
such that the  
homothety morphism
$$\chi^R_C  \colon R\to\rhom_R(C,C) $$
is an isomorphism is \emph{semidualizing}.  Observe that
the $R$-module $R$ is 
semidualizing.  
An  $R$-complex $D$ is \emph{dualizing} if it is 
semidualizing and has finite injective dimension;  
see Hartshorne~\cite[Chapter V]{hartshorne:rad} and
Foxby~\cite[Chapter 15]{foxby:hacr}.  
Over local rings,
dualizing complexes are unique up to shift-isomorphism.
\end{para}

The following result is proved like 
Jorgensen's~\cite[(2.5.1)]{jorgensen:jec}.

\begin{lem} \label{lem401}
Let $k$ be a field and $R_1,R_2$ local rings essentially of finite 
type over $k$ and let $R$ be a localization of $R_1\lotimes_k 
R_2$.  
If $D^i$ is a dualizing complex for $R_i$ for $i=1,2$, 
then the complex $(D^1\lotimes_k D^2)\lotimes_{R_1\lotimes_k 
R_2}R$ is dualizing for $R$.\qed
\end{lem}

\begin{para} \label{gd}
Let $C,X$ be homologically finite $R$-complexes with $C$ semidualizing.  
If the complex $\rhom_R(X,C)$ is homologically bounded and the 
biduality morphism 
$$\delta^C_X  \colon X\to\rhom_R(\rhom_R(X,C),C)$$
is an isomorphism, then $X$ is \emph{$C$-reflexive}.
The complexes $R$ 
and $C$ are $C$-reflexive,
and $C$ is dualizing if and only if each homologically finite 
complex is $C$-reflexive by~\cite[(V.2.1)]{hartshorne:rad}.
The 
\emph{$\text{G}_C$-dimension} of a $X$ 
is defined in~\cite{christensen:scatac} as
\[ \gkdim{C}_R(X)=
\begin{cases} \inf(C)-\inf(\rhom_R(X,C)) & \text{when $X$ is 
$C$-reflexive} \\
\infty & \text{otherwise.} \end{cases} \]
When $C=R$ this is the \emph{G-dimension} \label{para203}
of Auslander, Bridger, Foxby, and 
Yassemi~\cite{auslander:adgeteac,auslander:smt,yassemi:gd},
denoted 
$\gdim_R(X)$;  see also~\cite{christensen:gd}. 
If  $\pd_R(X)$ is finite, 
then
so is $\gdim_R(X)$,
and one has $\pd_R(\rhom_R(X,R))=-\inf(X)$ by~\cite[(2.13)]{christensen:scatac};
if in addition $R$ is local, then
$\gdim_R(X)=\pd_R(X)$ by~\cite[(2.3.10)]{christensen:gd}.
When $C,X$ are modules and $\gkdim{C}_R(X)=0$, one says $X$ is
\emph{totally $C$-reflexive}.
\end{para}

Other invariants and formulas are available over a local ring.

\begin{para} \label{para101}
When $R$ is local with residue field $k$ and $X$ is homologically 
finite, 
the integers
\begin{align*}
\beta_i^R(X)&=\rank_k(\HH_{-i}(\rhom_R(X,k)))
&\mu^i_R(X)&=\rank_k(\HH_{-i}(\rhom_R(k,X)) 
\end{align*}
are the $i$th \emph{Betti number} and \emph{Bass number} of 
$X$. 
The formal 
Laurent series
\begin{align*}
P^R_X(t)&=\sum_{i\in\zz}\beta_i^R(X) t^{i}
&I_R^X(t)&=\sum_{i\in\zz}\mu^i_R(X) t^i
\end{align*}
are the \emph{Poincar\'{e} series} and \emph{Bass series}
of $X$. The \emph{depth} of $X$ is
$$\depth_R(X)=-\sup(\rhom_R(k,X)).$$
When $C$ is a semidualizing $R$-complex, and $X$ is $C$-reflexive, the 
AB-formula reads
$$\gkdim{C}_R(X)=\depth(R)-\depth_R(X)$$
and the isomorphism
$R\simeq\rhom_R(C,C)$ gives rise to a formal equality
\[ P^R_C(t)I_R^C(t)=I^R_R(t)\]
by~\cite[(3.14)]{christensen:scatac} and~\cite[(1.5.3)]{avramov:rhafgd}.
When $D$ is dualizing for $R$, one has
$I_R^D(t)=t^d$ for some integer $d$ by~\cite[(V.3.4)]{hartshorne:rad}.
We say that $D$ is \emph{normalized} when $I_R^D(t)=1$, that is,
when $\inf(D)=\depth(R)$;  
see~\cite[(2.6)]{avramov:rhafgd}.  In particular, a minimal injective 
resolution $I$ of a normalized dualizing complex has $I_j\cong 
\oplus_{\p} E_R(R/\p)$ where the sum is taken over the set of prime 
ideals $\p$ with $\dim(R/\p)=j$.
\end{para}

We continue by recalling some standard morphisms.

\begin{para} \label{standard}
Let $X,Y,Z$ be $R$-complexes.  For an $R$-algebra $S$, let $U,V,W$ be 
$S$-complexes.
We have cancellation, commutativity, associativity, and adjunction 
isomorphisms.
\begin{align}
\label{sta} \tag{a}
X\lotimes_R R & \simeq X \\
\label{stb} \tag{b}
X\lotimes_R Y & \simeq Y\lotimes_R X \\
\label{stc} \tag{c}
X\lotimes_R (Y\lotimes_R Z) & \simeq (X\lotimes_R Y)\lotimes_R Z \\
\label{std} \tag{d}
\rhom_S(X\lotimes_R V,W) & \simeq \rhom_R(X,\rhom_S(V,W)) \\
\label{ste} \tag{e}
\rhom_R(U\lotimes_S V,Z) & \simeq \rhom_S(U,\rhom_R(V,Z))
\end{align}
Next, there are the tensor- and Hom-evaluation morphisms, 
respectively~\cite[(4.4)]{avramov:hdouc}.
\begin{align}
\label{stf} \tag{f}
\omega_{XVW}\colon 
& \rhom_R(X,V)\lotimes_S W \to \rhom_R(X,V\lotimes_S W) \\
\label{stg} \tag{g}
\theta_{XVW}\colon
& X\lotimes_R\rhom_S(V,W)\to \rhom_S(\rhom_R(X,V),W) 
\end{align}
The morphism $\omega_{XVW}$ is an isomorphism when
$X$ is homologically finite, $V$ is 
homologically bounded above, and either $\fd_S(W)<\infty$ or 
$\pd_R(X)<\infty$.  The morphism $\theta_{XVW}$ 
is an isomorphism when $X$ is homologically finite, 
$V$ is homologically bounded, and either $\injdim_S(W)<\infty$ or 
$\pd_R(X)<\infty$.
\end{para}

\begin{para} \label{para202}
Let 
$C,P,V,W,Y$ be $R$-complexes 
with
$Y$ homologically bounded above, 
$C$ semidualizing, 
and $\pd_R(P),\gkdim{C}_R(W)<\infty$.
\begin{enumerate}[\quad(a)]
\item \label{item216}
Adjunction and $C$-reflexivity provide an 
isomorphism 
\[ \rhom_R(V,W)\simeq\rhom_R(\rhom_R(W,C),\rhom_R(V,C)). \]
\item \label{item217}
Since $P$ is $R$-reflexive, 
Hom-evaluation 
gives an isomorphism
$$\rhom_R(P,Y)\simeq\rhom_R(P,R)\lotimes_RY.
$$
\end{enumerate}
\end{para}

In this paper we focus on several specific types of ring homomorphisms.

\begin{para}  \label{paraBass}
The ring homomorphism $\vf\colon R\to S$ induces a natural map on prime 
spectra $\vf^*\colon \spec(S)\to\spec(R)$.
The flat dimension of $\vf$ 
is defined as $\fd(\vf)=\fd_R(S)$.  

Assume 
that $\vf$ is \emph{local}, that is, the rings $R$ and $S$ 
are local with maximal ideals $\m$ and $\n$, respectively, and
$\vf(\m)\subseteq\n$.  The \emph{depth} of $\vf$
is $\depth(\vf)=\depth(S)-\depth(R)$.  When $\fd(\vf)$ is finite,
the \emph{Bass series} of $\vf$ is the formal Laurent 
series with nonnegative integer coefficients $I_{\vf}(t)$ satisfying 
the formal equality
$I_{S}^{S}(t)=I_R^R(t)I_{\vf}(t)$
whose existence is given by Avramov, Foxby, and 
Lescot~\cite[(5.1)]{avramov:bsolrhoffd} 
or~\cite[(7.1)]{avramov:rhafgd}.  The 
homomorphism $\vf$ is \emph{Gorenstein} at $\n$ if $I_{\vf}(t)=t^d$ for some 
integer $d$, in which case, $d=\depth(\vf)$.  When $\vf$ is 
module-finite, it is \emph{Cohen-Macaulay} if $S$ is perfect as an 
$R$-module, that is, when $\amp(\rhom_R(S,R))=0$.

\label{para01a}
When $\vf$ is surjective and has 
finite flat dimension (but is not necessarily local) 
it is
\emph{Cohen-Macaulay}
if, for each prime ideal $\q\subset S$, the 
localization $\vf_{\q}\colon R_{\p}\to S_{\q}$ is Cohen-Macaulay
where $\p=\vf^{*}(\q)$.
In this event, $\vf$ is \emph{Cohen-Macaulay of grade $d$}
if one of the following equivalent conditions holds:
\begin{enumerate}[\quad(i)]
\item
$S$ is a perfect $R$-module of grade $d$;
\item
$d=\grade_{R_{\p}}S_{\q}$ 
for each prime ideal $\q\subset S$;
\item
$\amp(\rhom_R(S,R))=0$.
\end{enumerate}
The map $\vf$ is \emph{Gorenstein}\footnote{Avramov
and Foxby~\cite{avramov:lgh,avramov:cmporh} 
originally used the terms
\emph{locally Cohen-Macaulay} and \emph{locally Gorenstein} for these
types of homomorphisms.  As they have chosen to rechristen the second
type \emph{Gorenstein}~\cite[(8.1)]{avramov:rhafgd}, 
we have followed suit with the first type.} 
if it is
Cohen-Macaulay and, 
for each prime ideal $\q\subset S$, the $S_{\q}$-module
$\ext^{d_{\q}}_R(S,R)_{\q}$ is cyclic for 
$d_{\q}=\grade_{R_{\p}}(S_{\q})$ where $\p=\vf^{*}(\q)$.
\end{para}

Here are two more combinations of standard morphisms.

\begin{para} \label{para402}
Assume that $\fd(\vf)$ is finite and fix $R$-complexes
$W,X,Y,Z$ 
with 
$W$ homologically bounded, $X$ homologically finite, and $Y$ 
homologically bounded above.
\begin{enumerate}[\quad(a)]
\item \label{item206}
Combining adjunction and tensor-evaluation yields an isomorphism
$$\rhom_S(X\lotimes_R S,Y\lotimes_R S)
\simeq\rhom_R(X,Y)\lotimes_R S.
$$
\item \label{item207}
If $\vf$ is module-finite, then
adjunction and Hom-evaluation provide
$$\rhom_S(\rhom_R(S,W),\rhom_R(S,Z))  
\simeq S\lotimes_R\rhom_R(W,Z). $$
\item \label{item208}
If $\vf$ is module-finite, then~\ref{para202}\eqref{item217},
tensor-evaluation, and
adjunction yield
$$\rhom_R(X\lotimes_RS,\rhom_R(S,Y))
\simeq \rhom_R(X,Y)\lotimes_R \rhom_R(S,R).$$
\end{enumerate}
\end{para}

When $X$ and $Y$ are modules the next lemma is~\cite[(2.5.8)]{grothendieck:ega4-2}.
Example~\ref{ex507} demonstrates the necessity of flatness.

\begin{lem} \label{sri}
Let $\vf\colon R\to S$ be flat and local such that the induced extension of residue
fields is bijective. 
If $X,Y$ are homologically degreewise finite and bounded below 
$R$-complexes and $X\lotimes_RS\simeq Y\lotimes_RS$ in $\D(S)$,
then $X\simeq Y$ in $\D(R)$.
\end{lem}

\begin{proof}
Consider minimal $R$-free resolutions $P\simeq X$ and $Q\simeq Y$.  
The $S$-complexes $P\otimes_RS$ and $Q\otimes_RS$
are minimal $S$-free resolutions of 
$X\lotimes_RS$ and $Y\lotimes_RS$, respectively.
The first
isomorphism in the following sequence follows from the flatness of $\vf$
\begin{align*}
\Hom_R(P,Q)\otimes_RS
&\simeq\rhom_R(X,Y)\lotimes_RS\\
&\simeq\rhom_{S}(X\lotimes_RS, Y\lotimes_RS)\\
&\simeq\Hom_{S}(P\otimes_RS, Q\otimes_RS)
\end{align*}
while the second is in~\ref{para402}\eqref{item206}
and the third is standard.
This shows that the composition of tensor-evaluation and adjunction
$$
f\colon \Hom_R(P,Q)\otimes_RS\xra{\simeq}\Hom_{S}(P\otimes_RS, Q\otimes_RS)
$$
is a quasiisomorphism.
The relevant definitions provide an equality
$$
\partial_0^{\Hom_R(P,Q)\otimes_RS}=\partial_0^{\Hom_R(P,Q)}\otimes_RS
$$
and the flatness of $\vf$ provides a natural isomorphism
\begin{equation} \label{eqker} \tag{$\dagger$}
\Ker(\partial_0^{\Hom_R(P,Q)})\otimes_R S
\cong\Ker(\partial_0^{\Hom_R(P,Q)\otimes_R S}).
\end{equation}
Note that the set of chain maps from $P$ to $Q$ over $R$ is exactly 
the set of  cycles  $Z_0^{\Hom_R(P,Q)}=\Ker(\partial_0^{\Hom_R(P,Q)})$, and
similarly for $P\otimes_RS$ and $Q\otimes_RS$.  

The assumption $X\lotimes_RS\simeq Y\lotimes_RS$
provides an isomorphism in the category of $S$-complexes
$\alpha\colon P\lotimes_RS\xra{\cong} Q\lotimes_RS$.
Since $f$ is a quasiisomorphism, there exists a cycle 
$\alpha'\in\Hom_R(P,Q)_0\otimes_RS$ such that 
the images of $f(\alpha')$ and $\alpha$ in $\HH_0(\Hom_{S}(P\otimes_RS, Q\otimes_RS))$
are equal.  In other words, the chain maps $f(\alpha')$ and $\alpha$ are homotopic.  
In particular, since $P\otimes_RS$ and $Q\otimes_RS$ are minimal and
$\alpha$ is an isomorphism of complexes, the same is true of $f(\alpha')$.

The isomorphism~\eqref{eqker} shows that $\alpha'=\sum_i \alpha_i'\otimes s_i$
for some $\alpha_i'\in\Ker(\partial_0^{\Hom_R(P,Q)})$ and $s_i\in S$.
For each $i$ fix an $r_i\in R$ with the same residue as $s_i$ in 
$k=R/\m\cong S/\n$.
We shall show that the chain map
$\alpha''=\sum_ir_i\alpha_i'\colon P\to Q$ is an isomorphism
of complexes.  
By construction, there is a commutative diagram
$$
\xymatrix{
(P\otimes_R S)\otimes_S k \ar[r]^{f(\alpha')\otimes_Sk}_{\cong} 
\ar[d]_{\cong}^{\ref{standard}\eqref{sta}\eqref{stc}} 
& (Q\otimes_R S)\otimes_S k \ar[d]_{\cong}^{\ref{standard}\eqref{sta}\eqref{stc}} \\
P\otimes_Rk \ar[r]^{\alpha''\otimes_Rk} & Q\otimes_Rk
}
$$
showing that $\alpha''\otimes_Rk$ is (degreewise) surjective.  
Nakayama's Lemma then implies that $\alpha''$ is degreewise
surjective, and the result follows from~\cite[(2.4)] {matsumura:crt}.
\end{proof}

The final background concept for this paper is the Picard group.

\begin{para} \label{para02}
The \emph{Picard group} of $R$, denoted $\Pic(R)$,
is the abelian group of
isomorphism classes of finitely generated locally free (i.e.,
projective) $R$-modules of rank 1 with operation given by tensor product. 
The assignment 
$M\mapsto M\lotimes_R S$
yields a well-defined group 
homomorphism $\Pic(\vf)\colon\Pic(R)\to\Pic(S)$.  
\end{para}

\section{Resolutions and localization} \label{sec5}

This section contains results used
to globalize standard local results.
We begin by observing that 
$\text{G}_C$-dimension can be measured by resolutions
when $C$ is a module. 
Example~\ref{ex02} shows that this fails when $\amp(C)>0$ 
however, see Lemma~\ref{lem209}.  

\begin{lem} \label{lem08}
Let $X$ be a homologically finite $R$-complex and $C$ a 
semidualizing $R$-module.
Given an integer $n$, the following conditions are equivalent:
\begin{enumerate}[\quad\rm(i)]
\item \label{item95c}
There is an isomorphism $G\simeq X$ where $G$ is a 
complex
of totally $C$-reflexive modules with $G_i=0$ for each $i>n$ and for
each $i<\inf (X)$;
\item \label{item95}
There is an inequality $\gkdim{C}_R(X)\leq n$;
\item \label{item95a}
One has $\gkdim{C}_R(X)<\infty$ and $n\geq-\inf(\rhom_R(X,C))$;
\item \label{item95b}
$n\geq\sup (X)$ and 
in any bounded below complex $G$ of totally $C$-reflexive modules
with $G\simeq X$,  
the module $\coker(\partial_{n+1}^G)$ is totally $C$-reflexive.
\end{enumerate}
In particular, there is an inequality $\sup(X)\leq\gkdim{C}_R(X)$.
\end{lem}

\begin{proof}
The local case when $X$ is a module 
is stated in~\cite[p.~68]{golod:gdagpi}.  For
the general case, mimic the proof of~\cite[(2.3.7)]{christensen:gd}.
\end{proof}

The next result is~\cite[(3.12)]{christensen:scatac} which we state here for ease of reference.
Example~\ref{ex501} shows that equality or strict inequality can occur.

\begin{lem} \label{lem209}
If $C,X$ are homologically finite $R$-complexes with 
$C$ semidualizing, then
$\sup (X)-\amp (C)\leq\gkdim{C}_R (X)$.\qed
\end{lem}

\begin{lem} \label{lem05}
If $C$ is homologically finite, the following conditions are 
equivalent:
\begin{enumerate}[\quad\rm(i)]
\item \label{item79a} $C$ is $R$-semidualizing; 
\item \label{item79b} $S^{-1}C$ is  $S^{-1}R$-semidualizing 
for each multiplicative subset $S\subset R$; 
\item \label{item79c} $C_{\m}$ is a $R_{\m}$-semidualizing 
for each maximal ideal $\m\subset R$.
\end{enumerate}
\end{lem}

\begin{proof}
The implication 
\eqref{item79b}$\implies$\eqref{item79c} is trivial, while
\eqref{item79a}$\implies$\eqref{item79b} follows from the argument 
of~\cite[(2.5)]{christensen:scatac}.   
For the remaining implication, condition~\eqref{item79c}
implies that the natural map 
$\chi^R_C\colon R\to\rhom_R(C,C)$ is locally an isomorphism,
so it is an isomorphism and $C$ is $R$-semidualizing.
\end{proof}

The proof of the next result
is almost identical to that 
of~\cite[(3.16)]{christensen:scatac}. 
Examples~\ref{ex502}--\ref{ex02} show that the inequalities can be strict or not,
that the converse of the second statement fails, and
that the final inequality fails to hold if $\amp(C)>0$.

\begin{lem} \label{lem05aa}
Let 
$C,X$ be homologically finite $R$-complexes with 
$C$ semidualizing.
For each multiplicative subset $S\subset R$, there is an inequality
$$\gkdim{S^{-1}C}_{S^{-1}R}(S^{-1}X)\leq 
\gkdim{C}_R(X)+\inf(S^{-1}C)-\inf(C).$$
In particular, if $X$ is $C$-reflexive, then $S^{-1}X$ is 
$S^{-1}C$-reflexive.
Furthermore, if $\amp(C)=0$, then  $\gkdim{S^{-1}C}_{S^{-1}R}(S^{-1}X)\leq 
\gkdim{C}_R(X)$.
\qed
\end{lem}

\begin{ex} \label{ex502}
When $R$ is local and $\amp(C)=0=\amp(X)$, the inequalities
in Lemmas~\ref{lem209}
and~\ref{lem05aa} can be strict (if $0\leq\pd_{S^{-1}R}(S^{-1}X)<\pd_R(X)<\infty$) or not
(set $C=R=X$).
\label{ex501}
\end{ex}

\begin{ex} \label{ex503}
The converse to the second statement in Lemma~\ref{lem05aa}
can fail.  Let $(R,\m)$ be a local non-Gorenstein ring
with prime ideal $\p\subsetneq \m$.  The module $\m$ is not $R$-reflexive
but the module $\m_{\p}\cong R_{\p}$ is 
$R_{\p}$-reflexive.  
\end{ex}

\begin{ex} \label{ex02}
The final inequalities in 
Lemmas~\ref{lem08} and~\ref{lem05aa} can fail if $\amp(C)>0$.  
Note that Lemma~\ref{lem101} shows that $X$ cannot be a semidualizing module.

Let $k$ be a field and $R=k[\![Y,Z]\!]/(Y^2,YZ)$.  Since $R$ is
complete local, it admits a dualizing complex $D$.  With $\p=(Y)R$ and
$X=R/\p$  the AB-formula implies 
\begin{gather*}
\gkdim{D}_R (X)=\depth (R)-\depth_R (X)=-1<0=\sup(X) 
\\
\gkdim{D_{\p}}_{R_{\p}} 
(X_{\p})=\depth(R_{\p})-\depth_{R_{\p}}(X_{\p})=0  
>-1=\gkdim{D}_R (X). 
\end{gather*}
The next equalities are by definition, and the first inequality is by~\cite[(A.8.6.1)]{christensen:gd}
\begin{gather*}
\gkdim{D}_R (D)=\inf (D)=\sup(D)-1<\sup(D) 
\\
\gkdim{D_{\p}}_{R_{\p}} (D_{\p})=\inf(D_{\p})=\inf(D)+1>\inf(D)  
=\gkdim{D}_R (D) 
\end{gather*}
while 
the second inequality follows from the arguments 
of~\cite[Section 15]{foxby:hacr}.  
\end{ex}

We do not know if the extra hypotheses are necessary
for the converses 
in the next result;
they are 
not needed when $\text{G}_C$-dimension is replaced by
projective dimension.

\begin{prop} \label{lem05ab}
Let 
$C,X$ be homologically finite $R$-complexes with 
$C$ semidualizing.
Consider the following conditions:  
\begin{enumerate}[\quad\rm(i)]
\item \label{item81a} $X$ is $C$-reflexive; 
\item \label{item81b} $S^{-1}X$ is $S^{-1}C$-reflexive
for each multiplicative subset $S\subset R$; 
\item \label{item81c} $X_{\m}$ is $C_{\m}$-reflexive
for each maximal ideal $\m\subset R$.
\end{enumerate}
The implications 
(i)$\implies$(ii)$\implies$(iii)
always 
hold, and 
the converses hold when either $\inf(\rhom_R(X,C))\geq-\infty$, 
$\dim(R)<\infty$, 
or $X$ is semidualizing.
\end{prop}

\begin{proof}
The implication 
\eqref{item81b}$\implies$\eqref{item81c} is trivial, while
\eqref{item81a}$\implies$\eqref{item81b} is in  
Lemma~\ref{lem05aa}.  
So, assume that 
$X_{\m}$ is $C_{\m}$-reflexive
for each maximal ideal $\m$. 
The biduality map
$\delta^C_X\colon X\to\rhom_R(\rhom_R(X,C),C)$ is locally an isomorphism,
and so it is an isomorphism.  It remains to show that
$\rhom_R(X,C)$ is homologically bounded. 

Assume first $\dim(R)<\infty$. 
For each maximal
ideal $\m\subset R$ the AB-formula provides the following equality 
while the inequality is~\cite[(2.7)]{foxby:daafuc}.
$$\gkdim{C_{\m}}_{R_{\m}}(X_{\m})=\depth 
(R_{\m})-\depth_{R_{\m}}(X_{\m})
\leq\dim (R_{\m})+\sup (X_{\m})$$
This explains the first inequality in the next sequence,
while the equality is by definition and the second inequality
is standard. 
\begin{align*}
\inf(\rhom_{R_{\m}}(X_{\m},C_{\m}))
& = \inf (C_{\m})-\gkdim{C_{\m}}_{R_{\m}}(X_{\m}) \\
& \geq \inf (C_{\m})-\dim (R_{\m})-\sup (X_{\m}) \\
& \geq \inf (C) - \dim (R) -\sup (X)
\end{align*}
It follows that $\rhom_R(X,C)$ is homologically bounded because
\begin{align*}
\inf(\rhom_R(X,C))
&=\inf\{\inf(\rhom_{R_{\m}}(X_{\m},C_{\m}))\mid\m\in\mspec(R)\}.
\end{align*}

Assuming next that $X$ is semidualizing,
the AB-formula and~\cite[(3.2.a)]{christensen:scatac}  provide
the equality 
$\gkdim{C_{\m}}_{R_{\m}}(X_{\m})
=\inf(X_{\m})$.
As above one deduces
$$\inf(\rhom_{R_{\m}}(X_{\m},C_{\m}))\geq \inf(C)-\sup(X)$$
and the homological boundedness of $\rhom_R(X,C)$.
\end{proof}

For strictness in the next inequality,
see
Example~\ref{ex106} or argue as in Example~\ref{ex102}.  

\begin{prop} \label{lem06b}
If $C$ is $R$-semidualizing, then there is an inequality
$$\gkdim{C}_R(X)\leq\sup\{\gkdim{C_{\m}}_{R_{\m}}(X_{\m})\mid\m\in\mspec(R)\}
$$
for each homologically finite $R$-complex $X$,
with equality if $\amp(C)=0$.
\end{prop}

\begin{proof}
For the inequality, set 
$s=\sup\{\gkdim{C_{\m}}_{R_{\m}}(X_{\m})\mid\m\in\mspec(R)\}$ and
 $i=\inf(\rhom_R(X,C))$, and
assume 
$s<\infty$.  For each maximal 
ideal $\m$, one has 
$$\gkdim{C_{\m}}_{R_{\m}}(X_{\m})+\inf(\rhom_R(X,C)_{\m})
=\inf(C_{\m})
\geq\inf(C). 
$$ 
It follows that $\rhom_R(X,C)$ is bounded because
the previous sequence gives 
$$i=\inf\{\inf(\rhom_R(X,C)_{\m})\mid 
\m\in\mspec(R)\}\geq \inf(C)-s
$$ 
so $\gkdim{C}_R(X)<\infty$ by 
Proposition~\ref{lem05ab}.
With 
$\m\in\supp_R(\HH_i(\rhom_R(X,C)))$, 
the desired inequality is in the next sequence.
$$\gkdim{C}_R(X)=\inf(C)-\inf(\rhom_R(X,C)) 
\leq\inf(C_{\m})-\inf(\rhom_R(X,C)_{\m})
\leq s$$
When $\amp(C)=0$, 
equality 
follows from Lemma~\ref{lem05aa} since 
$\inf(C)=\inf(C_{\m})$.
\end{proof}

Example~\ref{ex102}
shows the need for the connectedness hypothesis in the next result.
When $R$ is Cohen-Macaulay, the condition $\amp (C_{\m})=0$ is automatic
by~\cite[(3.4)]{christensen:scatac}.

\begin{prop} \label{prop101}
Let $C$ be a semidualizing $R$-complex and assume that $\spec(R)$
is connected.
If $\amp (C_{\m})=0$ for each maximal ideal $\m$,
then $\amp (C)=0$.
\end{prop}

\begin{proof}
If  $\amp(C)>0$, then $\spec(R)=\supp_R(C)$ is the disjoint union of the
closed sets $\supp_R(\HH_{\inf(C)}(C)),\ldots,\supp_R(\HH_{\sup(C)}(C))$,
contradicting
connectedness. 
\end{proof}

\begin{question}
If $C$ is a semidualizing $R$-complex and $\spec(R)$
is connected, must the inequality 
$\amp(C)=\sup\{\amp(C_{\m})\mid\m\in\mspec(R)\}$
hold?  
\end{question}

Proposition~\ref{prop101} 
with $C=\rhom_R(S,R)$ yields 
the next local-global principle; 
see~\ref{para01a}.  
Example~\ref{ex302} shows that this fails if $\spec(S)$ is disconnected. 

\begin{cor} \label{prop203}
Let $\vf\colon R\to S$ be a surjective Cohen-Macaulay
ring homomorphism. 
If $\spec(S)$ is connected, 
then 
$\vf$ is Cohen-Macaulay of constant grade.\qed
\end{cor}

\begin{ex} \label{ex302}
We construct a a surjective Cohen-Macaulay
ring homomorphism of nonconstant grade.
Let $k$ be a field and $R=k[Y,Z]$ a polynomial ring, and set
$$S=R/((Y,Z)R\cap (Y-1)R)$$ 
with
natural surjection $\vf\colon R\to S$.  Since $R$ is regular, 
one has $\pd_R(S)<\infty$. 
The equality $(Y,Z)R+(Y-1)R=R$ provides 
an 
isomorphism of $R$-algebras
$$S\cong R/(Y,Z)R \times R/(Y-1)R.$$
In particular, the ring $S$ is Cohen-Macaulay, and 
hence so is $\vf$
by~\cite[(8.10)]{avramov:cmporh}.  
Set $\n_1=(Y,Z)S$ and $\n_2=(Y-1,Z)S$.
To prove that $\vf$ has nonconstant grade, it suffices by~\ref{para01a} to show that 
$\amp(\rhom_R(S,R))>0$.  
For this
we verify 
\begin{align*}
\inf(\rhom_R(S,R)_{\n_1})&=-2
&\inf(\rhom_R(S,R)_{\n_2})&=-1.
\end{align*}
It is straightforward to verify that the localization $\vf_{\n_1}$ is 
equivalent to the natural surjection
$R_{(Y,Z)}\to k$ which has projective dimension 2.  Thus, one has
$$\inf(\rhom_R(S,R)_{\n_1})=\inf(\rhom_{R_{(Y,Z)}}(k,R_{(Y,Z)}))
=-\pd_{R_{(Y,Z)}}(k)=-2$$
where the second equality is by~\cite[(2.13)]{christensen:scatac}.
This is the first 
desired equality;  the second one follows similarly from the 
fact that the localization $\vf_{\n_2}$ is 
equivalent to the surjection
$k[Y,Z]_{(Y-1,Z)}\to k[Z]_{(Z)}$ which has projective dimension 1. 
\end{ex}

\begin{lem} \label{lem03}
Let $R=\coprod_{i\geq 0}R_i$ be a graded ring 
where $R_0$ is local with maximal ideal $\m_0$.  
Set $\m=\m_0+\coprod_{i\geq 1}R_i$ and let $X,Y$ 
be  homologically degreewise finite
complexes of graded $R$-module homomorphisms. 
\begin{enumerate}[\quad\rm(a)]
\item \label{item94}
For each integer $i$, one has $\HH_i(X)=0$ if and only if
$\HH_i(X_{\m})=0$.
\item \label{item94c}
There are equalities
$$ \inf(X)=\inf(X_{\m}) \quad \sup(X)=\sup(X_{\m})\quad
\amp(X)=\amp(X_{\m})$$
so  $X$ is homologically bounded (respectively,
homologically bounded above or homologically bounded below) 
if and only if the same is true of $X_{\m}$.
\item \label{item94d}
If $\alpha\colon X\to Y$ is a graded homomorphism of complexes, 
then $\alpha$ is a quasiisomorphism if and only if $\alpha_{\m}$
is a quasiisomorphism.
\end{enumerate}
\end{lem}

\begin{proof}
Part~\eqref{item94} 
follows from~\cite[(1.5.15)]{bruns:cmr}  
and the isomorphism $\HH_i(X_{\m})\cong\HH_i(X)_{\m}$,
and~\eqref{item94c} is immediate from~\eqref{item94}.
For~\eqref{item94d}, apply~\eqref{item94c} to the
mapping cone of $\alpha$.
\end{proof}

\begin{prop} \label{lem03a}
Let $R=\coprod_{i\geq 0}R_i$ be a graded ring 
where $R_0$ is local with maximal ideal $\m_0$.  
Set $\m=\m_0+\coprod_{i\geq 1}R_i$ and let $C,X$ 
be  homologically degreewise finite
complexes of graded $R$-module homomorphisms. 
\begin{enumerate}[\quad\rm(a)]
\item \label{item94b}
The complex $C$ is $R$-semidualizing if and only if $C_{\m}$
is $R_{\m}$-semidualizing.
\item \label{item94a}
If $C$ is $R$-semidualizing, then
$\gkdim{C}_R(X)=\gkdim{C_{\m}}(X_{\m})$.
Thus, the complex $X$ is $C$-reflexive if and only if
$X_{\m}$ is $C_{\m}$-reflexive.
\end{enumerate}
\end{prop}

\begin{proof}
\eqref{item94b}
One implication is contained in Lemma~\ref{lem05}, 
so assume that $C_{\m}$ is $R_{\m}$-semidualizing.
By Lemma~\ref{lem03}, the $R$-complexes $C$ 
and $\rhom_R(C,C)$ are homologically finite, 
and the homothety morphism
$R\to\rhom_R(C,C)$ is a quasiisomorphism.
so $C$ is semidualizing.

\eqref{item94a}  
It suffices to prove the final statement.  Indeed, if 
$X$ is $C$-reflexive and 
$X_{\m}$ is $C_{\m}$-reflexive, then the equality is a consequence of 
the following sequence
\begin{align*}
\gkdim{C}_R(X)
&=\inf(C)-\inf(\rhom_R(X,C))\\
&=\inf(C_{\m})-\inf(\rhom_{R_{\m}}(X_{\m},C_{\m}))\\
&=\gkdim{C_{\m}}(X_{\m})
\end{align*}
where the 
second equality is by Lemma~\ref{lem03}\eqref{item94c}, and the 
others are by definition.

For the final statement, one implication is in 
Lemma~\ref{lem05aa}, 
so assume that 
$X_{\m}$ is $C_{\m}$-reflexive.  
Lemma~\ref{lem03}
implies that $X$ and $\rhom_R(X,C)$ are homologically finite
and the the biduality 
morphism $X\to
\rhom_{R}(\rhom_{R}(X,C),C)$ is a
quasiisomorphism, so $X$ is $C$-reflexive.
\end{proof}

\section{Duality:  global results} \label{sec3}

This section is primarily devoted to 
reflexivity relations between semidualizing 
complexes in the nonlocal setting.  
We begin with 
a global version 
of~\cite[(3.1),(3.4)]{gerko:sdc}.

\begin{lem} \label{lem02}
Let 
$C,C'$ be semidualizing $R$-complexes.
\begin{enumerate}[\quad\rm(a)]
\item \label{item59}
If $C'$ is $C$-reflexive, then $\rhom_R(C',C)$ is semidualizing and
$C$-reflexive with
$\gkdim{C}(\rhom_R(C',C))=\inf(C)-\inf(C')$.
\item \label{item60}
If $C'$ 
is $C$-reflexive, then the evaluation morphism
$C'\lotimes_R \rhom_R(C',C)\to C$
is an isomorphism. 
\item \label{item61}
If $C\lotimes_R C'$ is semidualizing, then $C$ 
is $C\lotimes_R C'$-reflexive and
the evaluation morphism
$C\to\rhom_R(C',C\lotimes_R C')$
is an isomorphism.
\end{enumerate}
\end{lem}

\begin{proof}
Part~\eqref{item59} is contained in~\cite[(2.11)]{christensen:scatac}. 
For parts~\eqref{item60} and~\eqref{item61}, observe that
the maps are locally isomorphisms by~\cite[(3.1),(3.4)]{gerko:sdc} and are 
thus isomorphisms.
\end{proof}

The next result follows immediately from the local case; see~\cite[(5.3)]{takahashi:hiatsb}.
Example~\ref{ex504}
shows that $\rhom_R(C',C)\not\sim R$  in general;  
see also Example~\ref{ex102}.

\begin{lem} \label{prop202}
If $C,C'$ are  $R$-semidualizing, $C'$ is 
$C$-reflexive, and $C$ is $C'$-reflexive, then $\rhom_R(C',C)_{\m}\sim R_{\m}$
and $C'_{\m}\sim C_{\m}$
for each maximal 
ideal $\m$. \qed
\end{lem}

\begin{ex} \label{ex504}
One can have $\rhom_R(C',C)\not\sim R$ in Lemma~\ref{prop202}.
Assume that there exists
$[L]\in\Pic(R)$ with $[L]\neq[R]$.
If $C$ is a semidualizing $R$-complex, then so is $C'=C\lotimes_R L$.  Furthermore,
$C'$ is 
$C$-reflexive and $C$ is $C'$-reflexive.  However,
one has $\rhom_R(C,C')\simeq L\not\sim R$ and
$\rhom_R(C',C)\simeq \rhom_R(L,R)\not\sim R$.
\end{ex}

The next lemma follows directly from~\cite[(3.1),(3.2),(4.8.c)]{christensen:scatac}.
To see that the second and third inequalities can be strict and that the others
can be equalities, consult
Example~\ref{ex103} or argue as in Example~\ref{ex102}.  For strictness in 
the first and last inequalities, let $R$ be local and $\amp(C')>0$, and use
the guaranteed equalities.

\begin{lem} \label{lem101}
If $C,C'$ are $R$-semidualizing and $C'$ is $C$-reflexive, then
\begin{gather*}
\inf (C)-\sup (C')\leq\inf(\rhom_R(C',C)) \leq \inf (C)-\inf (C')\\
\inf (C')\leq \gkdim{C}_R(C')\leq\sup (C')
\end{gather*}
with equality in the second and third inequalities
if $R$ is local or $\amp(C')=0$.  In particular, if $C,C'$ are both 
modules, then $C'$ is totally $C$-reflexive. \qed
\end{lem}

Example~\ref{ex102}
shows the need for the connectedness hypothesis in the next result.

\begin{cor} \label{prop101a}
Let $C,C'$ be  semidualizing $R$-complexes with $\amp(C)=0$.
If $\spec(R)$
is connected
and $C'$ is $C$-reflexive,
then $\amp(C')=0$. 
\end{cor}

\begin{proof}
By Proposition~\ref{prop101}
we may assume 
that $R$ is local.  The first inequality in the next sequence
is in Lemma~\ref{lem08}
$$\sup (C')\leq\gkdim{C}_R(C')=\inf (C')\leq\sup (C')$$
while the equality is in~\cite[(3.1),(3.2)]{christensen:scatac}
and the last 
inequality is immediate. 
\end{proof}

\begin{question} \label{q502}
 If $\spec(R)$ is 
connected and $C,C'$ are semidualizing complexes such that $C'$ is 
$C$-reflexive, does the  inequality $\amp(C')\leq\amp(C)$ hold?
\end{question}

The answer is ``yes'' when $R$ is local and $C$ is 
dualizing for $R$ by~\cite[(3.4a)]{christensen:scatac}.
The next result resolves the 
local case when $C$ is not necessarily dualizing.
The inequality can be strict (e.g., if $\amp(C)>0$
and $C'=R$) 
or not (e.g., if $R$ is Cohen-Macaulay).
Consult 
Example~\ref{ex102} to see 
the need for connectedness.

\begin{cor} \label{cor301}
Let $R$ be local and $C,C'$ semidualizing $R$-complexes. If $C'$ 
is $C$-reflexive, then   $\amp(C')\leq\amp(C)$.
\end{cor}

\begin{proof}
Since $R$ is local, the equality in the following sequence is in
Lemma~\ref{lem101} 
$$\inf(C')=\gkdim{C}_R(C')\geq\sup(C')-\amp(C)$$
while the inequality is in Lemma~\ref{lem209}.
\end{proof}

\begin{ex} \label{ex102}
The conclusions of Proposition~\ref{prop101} and
Corollaries~\ref{prop101a} and~\ref{cor301} can
fail if $\spec(R)$ is not connected.
Let $k_1,k_2$ be fields and set $R=k_1\times k_2$.  
With $\m_1=0\times 
k_2$ and $\m_2=k_1\times 0$,
one has $\spec(R)=\{\m_1,\m_2\}$ and $R_{\m_i}\cong k_i\cong R/\m_i$
for $i=1,2$. Hence, $R$ is Gorenstein, and an $R$-complex is dualizing 
if and only if it is semidualizing.  
If $pq\neq 0$, the next equality and isomorphism are easily verified.
\begin{gather*}
\amp((\shift^ak_1^p)\times(\shift^bk_2^q))=|a-b| \\
\rhom_{R}((\shift^ak_1^p)\times(\shift^bk_2^q),
(\shift^ck_1^r)\times(\shift^dk_2^s)) \simeq
(\shift^{c-a}k_1^{pr})\times(\shift^{d-b}k_2^{qs})
\end{gather*}
It follows that 
$(\shift^ck_1^r)\times(\shift^dk_2^s)$ is dualizing if and 
only if $r=s=1$.  
So, the dualizing complex $C'=k_1\times\shift k_2$ is $R$-reflexive, and 
the next computations are routine
\begin{gather*}
C'_{\m_1}  \simeq k_1 
\qquad \qquad  C'_{\m_2}  \simeq 
\shift k_2 \notag \\
\amp(C'_{\m_i})=0<1=\amp(C')\\
\amp(R)=0<1=\amp(C')
\end{gather*}
\end{ex}

Here are the reflexivity relations between 
$\rhom_R(A,C)$ and $\rhom_R(B,C)$ 
when Lemma~\ref{lem02}\eqref{item59} guarantees that they
are semidualizing.  
Example~\ref{ex103} shows that the first inequality can be an equality and the second one
can be strict. 
The first one can also be strict: Use the guaranteed
equality when $R$ is local and $\amp(B)>0$.

\begin{prop} \label{propDD}
Let $A,B,C$ be
semidualizing $R$-complexes such that $A$ and $B$ are both $C$-reflexive
and set $(-)^{\dagger_C}=\rhom_R(-,C)$.
There are inequalities
\begin{align*}
\gkdim{A}_R(B) -\sup(B) 
&\leq
\gkdim{B^{\dagger_C}}_R(A^{\dagger_C})-\inf(C)+\inf(A)\\
&\leq\gkdim{A}_R(B) -\inf(B)  
\end{align*}
with equality at the second inequality when $R$ is local
or $\amp(B)=0$. 
In particular,
$B$ is $A$-reflexive if and only if 
$A^{\dagger_C}$ is $B^{\dagger_C}$-reflexive.
\end{prop}

\begin{proof}
It suffices to verify the final statement.  Indeed, if 
$B$ is $A$-reflexive and
$A^{\dagger_C}$ is $B^{\dagger_C}$-reflexive, then
Lemma~\ref{lem101} 
combined with~\ref{para202}\eqref{item216} provide
the desired inequalities and, when $R$ is local or
$\amp(B)=0$, 
the equalities.

Assume that $B$ is $A$-reflexive, and note
that $A^{\dagger_C}$ is homologically finite since $A$ is 
$C$-reflexive.
Employ 
the isomorphism from~\ref{para202}\eqref{item216} and the 
fact that $B$ is $C$-reflexive
to conclude that 
the complex $\rhom_R(A^{\dagger_C},B^{\dagger_C})$ is 
homologically 
bounded. 
Next, consider the 
following commutative diagram of morphisms of complexes.
\[ \xymatrix{
A^{\dagger_C} 
\ar[r]^-{\delta^{B^{\dagger_C}}_{A^{\dagger_C}}}\ar[d]^{\simeq}_{\ref{lem02}\eqref{item60}} 
& \rhom_R(\rhom_R(A^{\dagger_C},B^{\dagger_C}),B^{\dagger_C}) 
\ar[d]^{\simeq}_{\ref{para202}\eqref{item216}} \\
(B\lotimes_R\rhom_R(B,A))^{\dagger_C} & \rhom_R(B^{\dagger_C\dagger_C},
\rhom_R(A^{\dagger_C},B^{\dagger_C})^{\dagger_C}) 
\ar[d]^{\simeq} \\
\rhom_R(B,\rhom_R(B,A)^{\dagger_C}) \ar[u]_{\simeq}^{\ref{standard}\eqref{std}} 
& \rhom_R(B,\rhom_R(A^{\dagger_C},B^{\dagger_C})^{\dagger_C}) 
\ar[l]^{\simeq}_{\ref{para202}\eqref{item216}} 
}\]
The unmarked map
$\rhom_R(\delta^C_B,\rhom_R(A^{\dagger_C},B^{\dagger_C})^{\dagger_C})$
is an isomorphism since
$B$ is $C$-reflexive.
Thus, $\delta^{B^{\dagger_C}}_{A^{\dagger_C}}$ is an 
isomorphism and $A^{\dagger_C}$ is $B^{\dagger_C}$-reflexive.

The converse follows from the isomorphisms $A\simeq 
A^{\dagger_C\dagger_C}$ and $B\simeq 
B^{\dagger_C\dagger_C}$.
\end{proof}

\begin{ex} \label{ex103} 
Here we construct a ring $R$ with $\spec(R)$ connected 
demonstrating the following:  In Lemma~\ref{lem101} the first and fourth inequalities can 
be equalities and the other inequalities can be strict, and in Proposition~\ref{propDD}
the first inequality can 
be an equality and the other inequality can be strict.
Let $k$ be a field and set 
\begin{gather*}
A_1=k[X_1,Y_1]/(X_1^2,X_1Y_1)\qquad A_2=k[X_2,Y_2]/(X_2^2,X_2Y_2)\\
A=A_1\otimes_k A_2\cong k[X_1,Y_1,X_2,Y_2]/(X_1^2,X_1Y_1,X_2^2,X_2Y_2).
\end{gather*}
The natural maps
$\vf_i\colon A_i\to A$ are faithfully flat since they are obtained by applying
$-\otimes_k A_i$ to
the faithfully flat maps $k\to A_j$.  
For $i=1,2$ set $S_i=A_i\smallsetminus(X_i,Y_i)A_i$. The 
local ring $R_i=S_i^{-1}A_i$ has maximal ideal
$\m_i=(X_i,Y_i)R_i$ and exactly
one nonmaximal prime ideal $\p_i=(X_i)R_i$.  
Let $S=A\smallsetminus((X_1,Y_1,X_2)A\cup(X_1,X_2,Y_2)A)$ and set 
$R=S^{-1}A$ which has exactly two maximal ideals
$\n_1=(X_1,Y_1,X_2)R$ and $\n_2=(X_1,X_2,Y_2)R$ and exactly one 
nonmaximal prime ideal $\p=(X_1,X_2)R$.  As $\p\subset\n_1\cap\n_2$,
$\spec(R)$ is connected.

The containment $\vf_i(S_i)\subset S$ 
provides faithfully flat maps $\psi\colon R_i\to R$.  It is straightforward to verify that
$R$ is a localization of the tensor product $R_1\otimes_k R_2$,
and furthermore that
$\psi_i$ is the composition of the tensor product map
$R_i\to R_1\otimes_k R_2$ and the localization map
$R_1\otimes_k R_2\to R$.  Equally straightforward are the following.
\begin{align*}
\psi_1^{*}(\n_1)&=\m_1
&\psi_1^{*}(\n_2)&=\p_1
&\psi_1^{*}(\p)&=\p_1\\
\psi_2^{*}(\n_1)&=\p_2
&\psi_2^{*}(\n_2)&=\m_2
&\psi_2^{*}(\p)&=\p_1
\end{align*}
In particular, if $M_i$ is a nonzero $R_i$-module of finite length, then
the $R$-module $M_i\otimes_{R_i}R$ is nonzero with 
finite length because $\supp_R(M_i\otimes_{R_i}R)=\{\n_i\}$.

Since $R_i$ is essentially of finite type over $k$, it admits
a normalized dualizing complex $D^i$.
Hence, $\sup(D^i)=\dim(R_i)=1$ and $\inf(D^i)=\depth(R_i)=0$. 
From the structure of $\spec(R_i)$,
the minimal $R_i$-injective resolution of $D^i$ is of the form
$$D^i \simeq 0\to E_{R_{i}}(R_i/\p_i)\to 
E_{R_{i}}(R_i/\m_i)
\to 0.$$
In particular, the $R_i$-module
$\HH_0(D^i)$ has nonzero finite length.

Set $C^i=D^i\lotimes_{R_i}R$ which is semidualizing for $R$
by Theorem~\ref{lem01a1}. 
By flatness, we have
$\HH_j(C^i)\cong\HH_j(D^i)\otimes_{R_i}R$ for each integer $j$.
In particular, since the $R_i$-module
$\HH_0(D^i)$ has nonzero finite length, the 
$R$-module $\HH_0(C^i)$ has finite length and 
$\supp_R(\HH_0(C^i))=\{\n_i\}$.
Nakayama's lemma implies
$\HH_0(C^1)\otimes_R\HH_0(C^2)=0$.

Using Lemma~\ref{lem401} and the isomorphism
$C^1\lotimes_R C^2 \simeq(D^1\otimes_k D^2)\lotimes_{R_1\otimes_k R_2}R$
we conclude that $C^1\lotimes_R C^2$ is dualizing for $R$.
Write $D=C^1\lotimes_R C^2$.
In particular, $C^1,C^2$ are $D$-reflexive, and 
Lemma~\ref{lem02} provides isomorphisms
$$\rhom_R(C^1,D)\simeq C^2\qquad\qquad
\rhom_R(C^2,D)\simeq C^1.$$
We claim that $\inf (D)>0$.  Indeed,
since
$\inf (C^i)=0$ for $i=1,2$ one has
$$\inf (D)=\inf (C^1\lotimes_R C^2)>\inf (C^1)+\inf (C^2)= 0$$
where the inequality is~\cite[(A.4.15)]{christensen:gd} using
the last line of the previous paragraph.

We now show $\inf(D)=1=\amp(D)$.  
The K\"unneth formula $\HH(D^1\otimes_kD^2)=\HH(D^1)\otimes_k\HH(D^2)$ 
and the equalities $\sup(D^i)=1$
provide the next equality 
$$\sup(D)\leq\sup(D^1\otimes_kD^2)=2$$
while the inequality is due to the fact that $D$ is a localization of $D^1\otimes_kD^2$.
Since $\inf(D)\geq 1$, one has $0\leq\amp(D)\leq 1$, and so it suffices to verify $\amp(D)\geq 1$.
For this, note that the localizations $R_{\n_i}$ are not Cohen-Macaulay and therefore
one has $\amp(D)\geq\amp(D_{\n_i})\geq 1$.  
The desired computations now follow readily:
\begin{align*}
\inf(D)-\sup(C^1)&=\inf(\rhom_R(C^1,D)) <\inf (D)-\inf (C^1) \\
\inf (C^1)&<\gkdim{D}_R(C^1)=\sup(C^1)\\
\gkdim{D}_R(C^1)-\sup(C^1)&=
\gkdim{(C^1)^{\dagger_D}}_R(D^{\dagger_D})-\inf(D)+\inf(D)\\
&< \notag
\gkdim{D}_R(C^1) -\inf(C^1).
\end{align*}
\end{ex}

We next extend~\cite[(2.9)]{christensen:scatac}.
Example~\ref{ex106} shows that this inequality can be strict.

\begin{prop} \label{lem104}
If $C,X$ are homologically finite $R$-complexes 
with $C$ semidualizing, then there is an inequality
$$\gkdim{C}_R(X)\leq\pdim_R(X)$$
with equality when $\pdim_R(X)$ is finite
and either
$R$ is local or
$\amp(C)=0$.
\end{prop}

\begin{proof}
Assume 
that $\pdim_R(X)$ is finite.
The finiteness of $\gkdim{C}_R(X)$ is 
in~\cite[(2.9)]{christensen:scatac}, and the local case of the equality 
is~\cite[(3.5)]{christensen:scatac}. 
Proposition~\ref{lem06b} provides the inequality in the following sequence
\begin{align*}
\gkdim{C}_R(X)
&\leq\sup\{\gkdim{C_{\m}}_{R_{\m}}(X_{\m})\mid\m\in\mspec(R)\}\\
&=\sup\{\pd_{R_{\m}}(X_{\m})\mid\m\in\mspec(R)\}\\
&=\pd_R(X)
\end{align*}
while the first equality is by the local case and the second equality 
is classical.

Assume now that $\amp(C)=0$.  
Let $P\simeq X$ be
a projective resolution and set $g=\gkdim{C}_R(X)$.
Lemma~\ref{lem08} 
implies that  $G=\coker(\partial^P_{g+1})$ is totally $C$-reflexive,
and one checks locally (using the AB-formulas) that $G$ is projective.
\end{proof}

Next we 
extend~\cite[(3.17)]{christensen:scatac}.
Example~\ref{ex505} shows that the inequalities can be strict. 
Partial converses 
of the final statement and conditions guaranteeing equality are in 
Theorems~\ref{lem204} and~\ref{lem205}; 
to see that the converse can fail 
consult Example~\ref{ex301}.

\begin{prop} \label{lem203}
Let $C,P,X$ be homologically finite $R$-complexes 
with $C$ 
semidualizing and $\pd_R(P)$ finite.
There are inequalities
\begin{align*}
\gkdim{C}_R(X\lotimes_R P)\leq\gkdim{C}_R(X)+\pd_R(P)  \\
\gkdim{C}_R(\rhom_R(P,X))\leq\gkdim{C}_R(X)-\inf(P). 
\end{align*}
In particular, if $X$ is $C$-reflexive, then so are 
$X\lotimes_R P$ and $\rhom_R(P,X)$.
\end{prop}

\begin{proof}
The final statement is proved as 
in~\cite[(3.17)]{christensen:scatac}.  For the inequalities, assume 
that the complexes $X$,  $X\lotimes_R P$, and $\rhom_R(P,X)$
are $C$-reflexive.  Since $\pd_R(P)$ is finite, adjunction 
and~\ref{para202}\eqref{item217} 
yield an isomorphism  
\begin{equation*} 
\rhom_R(X\lotimes_R P,C)\simeq\rhom_R(P,R)\lotimes_R\rhom_R(X,C)
\end{equation*}
and so  the following sequence provides the first inequality.
\begin{align*}
\gkdim{C}_R(X\lotimes_R P)
&=\inf(C)-\inf(\rhom_R(X\lotimes_R P,C))\\
&=\inf(C)-\inf(\rhom_R(P,R)\lotimes_R\rhom_R(X,C))\\
&\leq\inf(C)-\inf(\rhom_R(P,R))-\inf(\rhom_R(X,C))\\
&=\gkdim{C}_R(X)+\pd_R(P)
\end{align*}
Similarly, 
the Hom-evaluation isomorphism
gives a sequence of (in)equalities
\begin{align*}
\gkdim{C}_R(\rhom_R(P,X))
&=\inf(C)-\inf(\rhom_R(\rhom_R(P,X),C))\\
&=\inf(C)-\inf(P\lotimes_R\rhom_R(X,C))\\
&\leq\inf(C)-\inf(P)-\inf(\rhom_R(X,C))\\
&=\gkdim{C}_R(X)-\inf(P)
\end{align*}
providing the second inequality. 
\end{proof}

\section{Ring homomorphisms of finite flat dimension:  Base change} \label{sec2}

In this section we study the interaction between the semidualizing and 
reflexivity properties and the functor $-\lotimes_RS$
where 
$\vf\colon R\to S$ is a ring homomorphism of finite flat dimension.
We begin with a more general situation~\cite[(A.4.15),(A.5.5)]{christensen:gd}
wherein the inequalities may be strict
(see Example~\ref{ex506})
or not (use $P=R$).

\begin{para} \label{para201}
If $X,P$ are $R$-complexes 
such that
$P\not\simeq 0$ is bounded and $\fd_R(P)$ is finite, then
\begin{align*}
\inf(X\lotimes_R P)&\geq\inf(X)+\inf(P)\\
\sup(X\lotimes_R P)&\leq\sup(X)+\fd_R(P)\\
\amp(X\lotimes_R P)&\leq\amp(X)+\fd_R(P)-\inf(P).
\end{align*}
\end{para}

Our nonlocal version of the 
amplitude inequality, based on~\cite{iversen:aifc} 
and~\cite[(3.1)]{foxby:daafuc}, is next.
It is Theorem I 
from the introduction and provides
inequalities complimentary to those in~\ref{para201}.
Example~\ref{ex506} shows that,
without the hypothesis on $\mspec(R)$, bounds of this ilk 
and the nontrivial ensuing implications need not hold, and that
the inequalities can be strict;
to see that they may not be strict, use $P=R$.

\begin{thm} \label{prop201} 
Let 
$P$ 
be a homologically 
finite $S$-complex 
with $\fd_R(P)$ finite and
such that $\vf^*(\supp_S(P))$ contains $\mspec(R)$.
For each homologically degreewise finite $R$-complex
$X$
there are inequalities
\begin{align*}
\inf(X\lotimes_R P)&\leq\inf(X)+\sup(P)\\
\sup(X\lotimes_R P)&\geq\sup(X)+\inf(P)\\
\amp(X\lotimes_R P)&\geq\amp(X)-\amp(P).
\end{align*}
In particular, 
\begin{enumerate}[\quad\rm(a)]
\item \label{item204}
$X\simeq 0$ if and only if $X\lotimes_R P\simeq 0$;
\item \label{item205}
$X$ is homologically bounded if and only if $X\lotimes_R P$ is so;
\item \label{item204b}
If $\amp(P)=0$, e.g., if $P=S$, then $\inf(X\lotimes_R 
P)=\inf(X)+\inf(P)$.
\end{enumerate}
\end{thm}

\begin{proof}
For the first inequality, it suffices to verify the following 
implication:  If $\HH_n(X)\neq 0$, then 
$\inf(X\lotimes_R P)\leq n+\sup(P)$.   Indeed, if $X\simeq 0$, 
that is, if $\inf(X)=\infty$, then 
the inequality is trivial.  If $\inf(X)$ is finite, then using 
$n=\inf(X)$ gives the desired inequality.  And if $\inf(X)=-\infty$, 
then taking the limit as $n\to-\infty$ gives the desired inequality.

Fix an integer $n$ and assume that $\HH_n(X)\neq 0$.  Thus, there 
is a maximal ideal $\m\in\supp_R(\HH_n(X))\subseteq\supp_R(X)$, and by 
assumption there exists a prime ideal $\p\in\supp_S(P)$ such that 
$\vf^{*}(\p)=\m$.  The local homomorphism $\vf_{\p}\colon R_{\m}\to 
S_{\p}$ and the complexes $X_{\m}$ and $P_{\p}$ satisfy the 
hypotheses of~\cite[(3.1)]{foxby:daafuc}, providing the second equality
in the following sequence wherein the inequalities are straightforward
$$\inf(X\lotimes_R P)
\leq\inf((X\lotimes_R P)_{\p})
=\inf(X_{\m}\lotimes_{R_{\m}}P_{\p})
=\inf(X_{\m})+\inf(P_{\p})
\leq n+\sup(P)$$
and the first equality follows from the isomorphism
$(X\lotimes_R P)_{\p}\simeq X_{\m}\lotimes_{R_{\m}}P_{\p}$.

The second inequality is verified similarly.  The third inequality 
is an immediate consequence of the first two, and 
statements~\eqref{item204}, \eqref{item205}, and~\eqref{item204b}
follow directly.
\end{proof}

The proof of the next result 
is nearly 
identical to that of~\cite[(2.10)]{iyengar:golh}, 
using $X=\cone(\alpha)$ in Theorem~\ref{prop201}.  
Example~\ref{ex506} shows that the extra hypotheses are necessary for the
nontrivial implication.  

\begin{cor} \label{para10}
Let
$P$ 
be a homologically 
finite $S$-complex
with $\fd_R(P)<\infty$ 
and $\mspec(R)\subseteq\vf^*(\supp_S(P))$. 
If
$\alpha$ is a morphism of homologically degreewise finite 
$R$-complexes, 
then $\alpha$ is an isomorphism if and only if 
$\alpha\lotimes_R P$ is 
so. \qed
\end{cor}

Here is 
a partial converse 
for Proposition~\ref{lem203}.
The first inequality can be strict:
use the guaranteed equality
with $R$ local and 
$\inf(P)<\pd_R(P)$.
Example~\ref{ex505} shows that the second 
inequality may be strict and the first one may not. 

\begin{thm} \label{lem204}
Let 
$C,P,X$ be homologically finite $R$-complexes with $C$ semidualizing,
$\pd_R(P)$ 
finite, and $\mspec(R)$ contained in $\supp_R(P)$.
There are inequalities
\begin{align*}
\gkdim{C}_R(X)+\inf(P)
&\leq\gkdim{C}_R(X\lotimes_R P)\\
&\leq\gkdim{C}_R(X)+\pd_R(P).  
\end{align*}
In particular, the complexes $X$ and
$X\lotimes_R P$ 
are $C$-reflexive simultaneously.  If 
$R$ is local 
or $\amp(\rhom_R(P,R))=0$, then the second inequality is an equality.
\end{thm}

\begin{proof}
First we verify that $X$ and
$X\lotimes_R P$
are $C$-reflexive simultaneously.  
Theorem~\ref{prop201}\eqref{item205}
and the isomorphism~\ref{para202}\eqref{item217}
imply
that the complexes
$\rhom_R(X,C)$ 
and
$\rhom_R(X\lotimes_R P,C)$ 
are homologically 
bounded simultaneously. 
The following commutative diagram from~\cite[(3.17)]{christensen:scatac}
$$ \xymatrix{
X\lotimes_R P \ar[r]^-{\delta^C_{X\lotimes_R P}} 
\ar[d]_{\delta^C_X\lotimes_R P} & \rhom_R(\rhom_R(X\lotimes_R P,C),C) 
\\
\rhom_R(\rhom_R(X,C),C)\lotimes_R P \ar[r]^-{\simeq}_{\ref{standard}\eqref{stg}} & 
\rhom_R(\rhom_R(P,\rhom_R(X,C)),C)\ar[u]_{\simeq}^-{\ref{standard}\eqref{std}}
}
$$
shows that $\delta^C_{X\lotimes_R P}$ and $\delta^C_X\lotimes_R P$ are 
isomorphisms simultaneously.  Corollary~\ref{para10} then 
implies that
$\delta^C_{X\lotimes_R P}$ and $\delta^C_X$ are 
isomorphisms simultaneously.

For the (in)equalities, we assume that $X$ and
$X\lotimes_R P$
are $C$-reflexive.  The first inequality 
is verified in the next 
sequence where (1) is by definition and~\ref{para202}\eqref{item217}
\begin{align*}
\gkdim{C}_R(X\lotimes_R P)
&\stackrel{(1)}{=}\inf(C)-\inf(\rhom_R(P,R)\lotimes_R\rhom_R(X,C))\\
&\stackrel{(2)}{\geq}\inf(C)-\sup(\rhom_R(P,R))-\inf(\rhom_R(X,C))\\
&\stackrel{(3)}{\geq}\gkdim{C}_R(X)-\pd_R(\rhom_R(P,R))\\
&\stackrel{(4)}{=}\gkdim{C}_R(X)+\inf(P)
\end{align*}
(2) is by 
Theorem~\ref{prop201}, (3) is 
standard, and (4)
is by~\ref{para203}.  
The second inequality is in Proposition~\ref{lem203}.  
When $\amp(\rhom_R(P,R))=0$, 
there is an equality
$$\inf(\rhom_R(P,R)\lotimes_R\rhom_R(X,C))=
\inf(\rhom_R(P,R))+\inf(\rhom_R(X,C))$$
by Theorem~\ref{prop201}\eqref{item204b};  the same equality 
holds by Nakayama's Lemma when $R$ is local.  Thus, under either of 
these hypotheses, the displayed sequence in the proof of 
Proposition~\ref{lem203} gives the desired equality.
\end{proof}

The next result 
contains Theorem II\eqref{itemIa} from the introduction.
Example~\ref{ex506} shows that the converse of the first implication 
can fail.
To see that the inequalities can be strict, 
consult Example~\ref{ex106}.
For equality, use $C=R$.

\begin{thm} \label{lem01a1}
Assume 
that $\fd(\vf)$ is finite and $C$ is a  
homologically degreewise finite $R$-complex. 
If $C$ is $R$-semidualizing, then 
$C\lotimes_R S$ is $S$-semidualizing with
$$\inf(C\lotimes_R S)\geq\inf(C) \qquad \text{and} \qquad 
\amp(C\lotimes_R S)\leq\amp(C).$$
Conversely, if $C\lotimes_R S$ is 
$S$-semidualizing
and $\image(\vf^*)$ contains $\mspec(R)$,
then $C$ is $R$-semidualizing and 
$\inf(C\lotimes_R S)=\inf(C)$.
\end{thm}

\begin{proof}
The first implication and the inequalities are
in~\cite[(1.3.4),(5.1)]{christensen:scatac}.
Assume that $C\lotimes_R S$ is $S$-semidualizing
and $\mspec(R)\subseteq\image(\vf^*)$.
The equality is in 
Theorem~\ref{prop201}\eqref{item204b},
and Theorem~\ref{prop201}\eqref{item205} implies that $C$ is homologically 
finite over $R$.
The following commutative diagram shows that 
$\chi_C^R\lotimes_R S$ is an isomorphism
\[ \xymatrix{
S \ar[r]^-{\chi_{C\lotimes_R S}^S}_-{\simeq}\ar[d]^{\simeq}_{\ref{standard}\eqref{sta}} 
 & \rhom_S(C\lotimes_R S,C\lotimes_R S) \ar[d]^{\simeq}_{\ref{para402}\eqref{item206}} \\
R\lotimes_R S \ar[r]^-{\chi_{C}^R\lotimes_R S} & \rhom_R(C,C)\lotimes_R S
} \]
and Corollary~\ref{para10} 
implies that $\chi_C^R$ 
is an isomorphism as well.  
\end{proof}

\begin{ex} \label{ex506}
We show:
(1) the implication in Theorem~\ref{prop201}\eqref{item204b} can fail
when $\amp(P)>0$;
(2) the nontrivial implications in Theorem~\ref{prop201}\eqref{item204} 
and Corollary~\ref{para10} can fail in the absence 
of the hypothesis on $\supp_S(P)$; and
(3) the first inequality in Proposition~\ref{lem203} and the
second inequality in Theorem~\ref{lem204} can be strict, while equality can occur
in the first inequality
in Theorem~\ref{lem204}.

Let $R=k[Y]$.
Setting  $P^1=R/(Y-1)\oplus \shift R \oplus \shift^2 R/Y$ and $X^1=R/(Y)$,
one has $X^1\lotimes_R P^1\simeq \shift R/(Y)$, 
and so one verifies (1) from the next computations.
\begin{gather*}
\inf(X^1)=\sup(X^1)=\amp(X^1)=0 \\
\inf(P^1)=0 \qquad \sup(P^1)=\amp(P^1)=2\qquad \fd_R(P^1)=3 \\
\inf(X^1\lotimes_R P^1)=\sup(X^1\lotimes_R P^1)=1\qquad\amp(X^1\lotimes_R P^1)=0
\end{gather*}
For (2)
let $\alpha\colon R\to R/(Y)\oplus R$ be the natural map and $P^2=R/(Y-1)$.  
It is straightforward to check that $\alpha\lotimes_RP^2$ is an isomorphism,
even though $\alpha$ is not. 
Furthermore, with $X^2=\cone(\alpha)$ one has 
$X^2\lotimes_R P^2\simeq 0$ while
$X^2\not\simeq 0$.
With $X^3=R\oplus(\oplus_{i\in\mathbb{Z}}R/(Y))$, the complex $X^3\lotimes_R P^2\simeq P^2$
is homologically bounded, even though $X^3$ is not.

\label{ex505}
For (3), if $P^3= R\oplus R/(Y-1)$, then 
$X^1\lotimes_R P^3\simeq R/(Y)$ and so
$$\gkdim{C}_R(X^1)+\inf(P^3)=\gkdim{C}_R(X^1\lotimes_R P^3)
<\gkdim{C}_R(X^1)+\pdim_R(P^3).$$
Set $S=R/(Y)$ with $\vf\colon R\to S$ the 
natural surjection.  The module $P^3$ is not 
$R$-semidualizing, even though  $P^3\lotimes_R S\simeq S$ 
is $S$-semidualizing.  
\end{ex}

Next we refine the ascent property~\cite[(5.10)]{christensen:scatac}. 
When $\vf$ is local, Theorem~\ref{lem01a3} shows that this inequality can be strict
(if $\amp(C)>0$) or not (if $\amp(C)=0$).
Example~\ref{ex301} shows that the converse to the final 
statement need not hold.

\begin{prop} \label{lem01a2}
Assume that $\fd(\vf)$ is finite, and
let $C,X$ be homologically finite $R$-complexes
such that $C$ is $R$-semidualizing.
There is an inequality
$$\gkdim{C\lotimes_R S}_S(X\lotimes_R S)\leq\amp(C)+\gkdim{C}_R(X).$$
In particular, if $X$ is $C$-reflexive, then 
$X\lotimes_R S$
is $C\lotimes_R S$-reflexive.  
\end{prop}

\begin{proof}
The last statement is in~\cite[(5.10)]{christensen:scatac}, so assume
that $\gkdim{C}_R(X)$ and 
$\gkdim{C\lotimes_R S}_S(X\lotimes_R S)$ are finite.  
In the following sequence
\begin{align*}
\gkdim{C\lotimes_R S}_S(X\lotimes_R S)
&=\inf(C\lotimes_R S)-\inf(\rhom_S(X\lotimes_R S,C\lotimes_R S))\\
&\leq\sup(C)-\inf(\rhom_R(X,C))\\
&=\amp(C)+\inf(C)-\inf(\rhom_R(X,C))\\
&=\amp(C)+\gkdim{C}_R(X)
\end{align*}
the equalities are routine, and
the inequality follows from \ref{para402}\eqref{item206}
and~\ref{para201}.
\end{proof}
 
The following descent result is
Theorem II\eqref{itemIb} from the introduction.

\begin{thm} \label{lem01a3}
Let $C,X$ be homologically degreewise finite $R$-complexes
with $C$ semidualizing.
If $\fd(\vf)$ is finite and $\image(\vf^*)$ contains 
$\mspec(R)$,
then
$$\gkdim{C}_R(X)=\gkdim{C\lotimes_R S}_S(X\lotimes_R S).
$$ 
In particular, 
$X\lotimes_R S$
is $C\lotimes_R S$-reflexive if and only if $X$ is $C$-reflexive.
\end{thm}

\begin{proof}
One  implication is
in Proposition~\ref{lem01a2}, so 
assume that $X\lotimes_R S$
is $C\lotimes_R S$-reflexive.
Theorem~\ref{prop201}\eqref{item205}
and~\ref{para402}\eqref{item206}
imply that $X$ and $\rhom_R(X,C)$ are 
homologically bounded.  
With Corollary~\ref{para10} the commutative diagram
from~\cite[(5.10)]{christensen:scatac}
\[ \xymatrix{
X\lotimes_R S \ar[r]^-{\delta^{C\lotimes_R S}_{X\lotimes_R S}}_-{\simeq}
\ar_=[dd] &
\rhom_S(\rhom_S(X\lotimes_R S,C\lotimes_R S),C\lotimes_R S) \\
 & \rhom_S(\rhom_R(X,C)\lotimes_R S,C\lotimes_R S) \ar[u]_{\simeq}^{\ref{para402}\eqref{item206}} \\
X\lotimes_R S \ar[r]^-{\delta^C_{X}\lotimes_R S} 
& \rhom_R(\rhom_R(X,C),C)\lotimes_R S \ar[u]_{\simeq}^{\ref{para402}\eqref{item206}}
} \]
shows that
$\delta^C_{X}$ is an isomorphism,
and  so $X$ is $C$-reflexive.

Assuming that 
$\gkdim{C}_R(X)$ 
and
$\gkdim{C\lotimes_R S}_S(X\lotimes_R S)$ 
are finite, one has  
\begin{align*}
\gkdim{C\lotimes_R S}_S(X\lotimes_R S)
&=\inf(C\lotimes_R S)-\inf(\rhom_S(X\lotimes_R S,C\lotimes_R S)) \\
&=\inf(C)-\inf(\rhom_R(X,C))\\
&=\gkdim{C}_R(X)
\end{align*}
where the second equality
is from Theorem~\ref{prop201}\eqref{item204b} and~\ref{para402}\eqref{item206}. 
\end{proof}

Here is Theorem II\eqref{itemIc} from the introduction.
It uses the functor $\Pic(-)$; see~\ref{para02}.
The conclusion fails outright if $C,C'$ are not 
semidualizing by Example~\ref{ex507}.
Note that the injectivity of $\Pic(\vf)$ in the hypotheses is not automatic, 
even when $\vf$ is faithfully flat~\cite[(11.8)]{fossum:dcgkd}, unless $\vf$ 
is local or surjective; see Proposition~\ref{lem09}.
In fact, the 
inclusion $\Pic(R)\subseteq\s(R)$ shows that this condition is necessary.

\begin{thm} \label{lem01a4}
Assume that $\fd(\vf)$ is finite, $\image(\vf^*)\supseteq\mspec(R)$,
and $\Pic(\vf)$ is injective.
If $C,C'$ are $R$-semidualizing 
and $C\lotimes_R S\simeq C'\lotimes_R S$, then $C\simeq C'$.
\end{thm}

\begin{proof}
By Theorem~\ref{lem01a3} the isomorphism $C\lotimes_R S\simeq C'\lotimes_R S$ 
implies that $C$ is $C'$-reflexive and vice 
versa.  It follows from Lemma~\ref{prop202} that, for each 
$\m\in\mspec(R)$, there is an isomorphism $\rhom_R(C',C)_{\m}\sim 
R_{\m}$.  The isomorphisms
\begin{equation} \label{eq15} \tag{$\dagger$}
R\lotimes_R S\simeq S\simeq\rhom_S(C'\lotimes_R S,C\lotimes_R S)
\simeq\rhom_R(C',C)\lotimes_R S \end{equation}
along with
Theorem~\ref{prop201} explain the following inequalities
$$0=\amp(\rhom_R(C',C)\lotimes_R S)\geq\amp(\rhom_R(C',C))\geq 0.$$
Thus, $\amp(\rhom_R(C',C))=0$ and 
$\rhom_R(C',C)_{\m}\simeq 
\shift^iR_{\m}$ for each 
$\m\in\mspec(R)$, where $i=\inf(\rhom_R(C',C))$.  
In other words, $\rhom_R(C',C)\simeq\shift^i L$ 
where 
$[L]\in\Pic(R)$. The isomorphisms~\eqref{eq15} 
imply
$S\simeq\shift^iL\lotimes_R S\simeq\shift^iL\otimes_R S$
and so $i=0$.
Applying~\eqref{eq15} again yields
$\Pic(\vf)([L])=[S]=\Pic(\vf)([R])$
so the injectivity of $\Pic(\vf)$ implies $L\cong R$.  
Hence, $\rhom_R(C',C)\simeq R$ and thus
$$C'\simeq R\lotimes_R C'\simeq\rhom_R(C',C)\lotimes_R C'\simeq C$$
where the last isomorphism is from
Lemma~\ref{lem02}\eqref{item60}.
\end{proof}

\begin{ex} \label{ex507}
The conclusions of Lemma~\ref{sri} and Theorem~\ref{lem01a4} fail
if $\vf$ is not flat and if the complexes are not semidualizing.
Set $R=k[\![Y,Z]\!]$ and $S=R/(Y,Z)$ with $\vf\colon R\to S$ 
the surjection.  The 
complexes $C=R/(Y)$ and $C'=R/(Z)$ satisfy 
$C\lotimes_R S \simeq S\oplus\shift S\simeq C'\lotimes_R S$
and $C\not\sim C'$.
\end{ex}

\begin{prop} \label{lem09}
If $\vf$ is surjective with $\fd(\vf)$ finite and 
$\mspec(R)$ 
is 
contained in $\image(\vf^*)$, 
then $\Pic(\vf)$ is injective.
\end{prop}

\begin{proof}
Set $I=\Ker(\vf)$ so that $S\cong R/I$, and note that our hypothesis
on $\vf$ implies that the Jacobson radical of $R$ contains $I$.
Let $L$ be a finitely generated rank 1 projective $R$-module such that
$S\cong L\otimes_R S\cong L/IL$.  Fix an element $x\in L$ whose residue
in $L/IL$ is a generator and let $\alpha\colon R\to L$ be given by
$1\mapsto x$.  By construction, the induced map 
$\alpha\otimes_R S\colon S\to L\otimes_R S$ is bijective.
Since $L$ is a projective $R$-module, this says that the morphism
$\alpha\lotimes_R S\colon S\to L\lotimes_R S$ is an isomorphism.  
By Corollary~\ref{para10} it follows that $\alpha$ is also an 
isomorphism.
\end{proof}

\section{Finite ring homomorphisms of finite flat dimension:  Cobase 
change}
\label{sec7}

Here we study the relation between 
the semidualizing and 
reflexivity properties 
and the functor
$\rhom_R(S,-)$ where 
$\vf\colon R\to S$ is a module-finite
ring homomorphism of finite flat dimension.
We begin with results 
that follow directly 
from~\ref{para201}--\ref{para10}
using~\ref{para203} and~\ref{para202}\eqref{item217}; their 
limitations are shown by the same examples.

\begin{para} \label{para201a}
If 
$X,P$ are $R$-complexes 
such that $\HH(P)\neq 0$ is finite and $\pd_R(P)<\infty$, then
\begin{align*}
\inf(\rhom_R(P,X))&\geq\inf(X)-\pd_R(P)\\
\sup(\rhom_R(P,X))&\leq\sup(X)-\inf(P)\\
\amp(\rhom_R(P,X))&\leq\amp(X)+\pd_R(P)-\inf(P).
\end{align*}
\end{para}

\begin{cor} \label{prop201a}  
Let 
$P$ 
be a homologically 
finite $R$-complex 
with $\pd_R(P)$ finite and
such that $\supp_R(P)$ contains 
$\mspec(R)$. 
For each homologically degreewise finite $R$-complex
$X$
there are inequalities
\begin{align*}
\inf(\rhom_R(P,X))&\leq\inf(X)+\sup(\rhom_R(P,R))\\
\sup(\rhom_R(P,X))&\geq\sup(X)+\inf(\rhom_R(P,R))\\
\amp(\rhom_R(P,X))&\geq\amp(X)-\amp(\rhom_R(P,R)).
\end{align*}
In particular, 
\begin{enumerate}[\quad\rm(a)]
\item \label{item204a}
$X\simeq 0$ if and only if $\rhom_R(P,X)\simeq 0$;
\item \label{item205a}
$X$ is homologically bounded if and only if $\rhom_R(P,X)$ is so;
\item \label{item205b}
If $\amp((\rhom_R(P,R))=0$, then 
the first inequality is an equality.\qed
\end{enumerate}
\end{cor}

\begin{cor} \label{para10a}
Let 
$P$ 
be a homologically 
finite $R$-complex 
with $\pd_R(P)$ finite and such that
$\supp_R(P)$ contains $\mspec(R)$.
If
$\alpha$ is a morphism of homologically degreewise finite 
$R$-complexes, 
then $\alpha$ is an isomorphism if and only if the 
induced morphism $\rhom_R(P,\alpha)$ is an isomorphism.\qed
\end{cor}

Here is 
a partial converse 
for Proposition~\ref{lem203}.
For strictness in the first inequality, 
use the guaranteed equality
with $R$ local and 
$\amp(P)>0$.
Example~\ref{ex505} shows other limitations.

\begin{thm} \label{lem205}
Let 
$C,P,X$ be homologically finite $R$-complexes with $C$ semidualizing,
$\pd_R(P)$ 
finite, and $\mspec(R)$ contained in $\supp_R(P)$.
There are inequalities
\begin{align*}
\gkdim{C}_R(X)-\sup(P) 
&\leq\gkdim{C}_R(\rhom_R(P,X))\\
&\leq\gkdim{C}_R(X)-\inf(P). 
\end{align*}
In particular, the complexes $X$ and $\rhom_R(P,X)$
are $C$-reflexive simultaneously.  If 
$R$ is local 
or $\amp(P)=0$ then the second inequality is an equality.
\end{thm}

\begin{proof}
Set $(-)^{\dagger_C}=\rhom_R(-,C)$.
First we verify that $X$ and
$\rhom_R(P,X)$
are $C$-reflexive simultaneously. 
Theorem~\ref{prop201}\eqref{item205} and the Hom-evaluation isomorphism 
\begin{equation} \label{eq208} \tag{$\dagger$}
\rhom_R(P,X)^{\dagger_C}\simeq P\lotimes_RX^{\dagger_C}
\end{equation}
show that the complexes
$\rhom_R(P,X)^{\dagger_C}$ and $X^{\dagger_C}$
are homologically 
bounded simultaneously.
The following commutative diagram from~\cite[(3.17)]{christensen:scatac}
$$ \xymatrix{
\rhom_R(P,X) \ar[rr]^-{\delta^C_{\rhom_R(P,X)}} 
\ar[d]_{\rhom_R(P,\delta^C_X)} && \rhom_R(P,X)^{\dagger_C\dagger_C}
\ar[d]_{\simeq}^{\ref{standard}\eqref{stg}} \\
\rhom_R(P,X^{\dagger_C\dagger_C}) \ar[rr]^-{\simeq}_-{\ref{standard}\eqref{std}} && 
(P\lotimes_RX^{\dagger_C})^{\dagger_C}
}
$$
implies that $\delta^C_{\rhom_R(P,X)}$ and $\rhom_R(P,\delta^K_X)$ are 
isomorphisms simultaneously.  Corollary~\ref{para10a} shows 
that
the same is true for
$\delta^C_{\rhom_R(P,X)}$ and $\delta^C_X$.

Assume that $X$ and
$\rhom_R(P,X)$
are $C$-reflexive.  
The second inequality is in Proposition~\ref{lem203}.  
The first inequality 
is verified in the following 
sequence
\begin{align*}
\gkdim{C}_R(\rhom_R(P,X))
&\stackrel{(1)}{=}\inf(C)-\inf(P\lotimes_R\rhom_R(X,C))\\
&\stackrel{(2)}{\geq}\inf(C)-\sup(P)-\inf(\rhom_R(X,C))\\
& 
=\gkdim{C}_R(X)-\sup(P)
\end{align*}
where (1) is by 
isomorphism~\eqref{eq208}, and (2) is by 
Theorem~\ref{prop201}. 
If $\amp(P)=0$, then
$$\inf(P\lotimes_R\rhom_R(X,C))=
\inf(\rhom_R(P,R))+\inf(\rhom_R(X,C))$$
by Theorem~\ref{prop201}\eqref{item204b};  the same equality 
holds by Nakayama's Lemma if $R$ is local.  Thus, under either of 
these hypotheses, the displayed sequence in the proof of 
Proposition~\ref{lem203} gives the desired equality.
\end{proof}

Example~\ref{ex506} shows how the converse of the first implication 
of the next result can fail. 
If $R$ is local and $\amp(C)=0$, then the second and third inequalities are strict
if and only if $\pd_R(S)>0$.  We do not know if the first inequality can be strict.

\begin{thm} \label{lem102}
Assume that $\vf$ is module-finite with $\fd(\vf)$ finite
and $C$ is a 
homologically degreewise finite $R$-complex. 
If $C$ is $R$-semidualizing, then $\rhom_R(S,C)$ is 
$S$-semidualizing and 
$$\inf(C)-\pd_R(S)\leq\inf(\rhom_R(S,C))\leq\sup (C)$$
with equality on the left if $R$ is local or $\amp(C)=0$.
Conversely, if $\rhom_R(S,C)$ is 
$S$-semidualizing and $\mspec(R)\subseteq\image(\vf^*)$, 
then $C$ is $R$-semidualizing and 
$$\inf(\rhom_R(S,C))\leq\inf(C)+\sup(\rhom_R(S,R))$$ 
with equality 
if $\amp(\rhom_R(S,R))=0$.
\end{thm}

\begin{proof}
First, assume that $C$ is $R$-semidualizing.  
Mimic the proof 
of~\cite[(6.1)]{christensen:scatac}
to show that
$\rhom_R(S,C)$ is 
$S$-semidualizing.
The first inequality and conditional equality follow immediately from
Proposition~\ref{lem104}.
The second inequality is a consequence of~\ref{para201a} since
$\inf(\rhom_R(S,C))\leq\sup(\rhom_R(S,C))$. 

Next, assume that $\rhom_R(S,C)$ is $S$-semidualizing and
$\mspec(R)\subseteq\image(\vf^*)$.
Corollary~\ref{prop201a}\eqref{item205a} implies that $C$ 
is homologically finite.
The commutative diagram 
$$
\xymatrix{
S \ar[rr]^-{\chi^S_{\rhom_R(S,C)}}_-{\simeq} \ar[d]_{\simeq} &&
\rhom_S(\rhom_R(S,C),\rhom_R(S,C)) \ar[d]_{\simeq}^{\ref{para402}\eqref{item207}} \\
S\lotimes_R R 
\ar[rr]^-{S\lotimes_R\chi^R_C} 
&& S\lotimes_R\rhom_R(C,C)
}
$$
shows that 
$S\lotimes_R\chi^R_C$ is an isomorphism and  Corollary~\ref{para10} implies the same for
$\chi^R_C$. 
The last (in)equality is in Corollary~\ref{prop201a}.
\end{proof}

Part~\eqref{item38} of 
the next result says,  
if $\amp(\rhom_R(S,R))=0=\amp(C)$, then $\amp(\rhom_R(S,C))=0$. 

\begin{prop} \label{prop01a}
Let  $C$ be a semidualizing $R$-module, and assume that 
$\vf$ is surjective and Cohen-Macaulay of grade $d$.
\begin{enumerate}[\quad\rm(a)]
\item \label{item38}
$\ext^d_R(S,C)$ is 
$S$-semidualizing and $\ext^i_R(S,C)=0$ for each $i\neq d$.
\item \label{item32}
If $\vf$ is 
Gorenstein, then the $S$-module $\ext^d_R(S,R)$ is 
locally free of rank 1.
\end{enumerate}
\end{prop}

\begin{proof}
\eqref{item38}  
Let  $\q\subset S$ be prime and
set $\p=\vf^{*}(\q)$ and $I=\Ker(\vf)$.  
The $S$-complex 
$\rhom_R(S,C)$ is semidualizing by Theorem~\ref{lem102},
so it suffices to show $\ext_R^j(S,C)_{\q}=0$ for $j\neq d$.
There is an $R_{\p}$-sequence 
$\mathbf{y}\in I_{\p}$ of length $d=\grade_{R_{\p}}(S_{\q})$.
Since $C_{\p}$ is $R_{\p}$-semidualizing, 
$\mathbf{y}$
is also $C_{\p}$-regular, and thus 
$\ext_R^j(S,C)_{\q}=\ext_{R_{\p}}^j(S_{\q},C_{\p})=0$
for $j<d$.
Also, $d=\pdim_{R_{\p}}(S_{\q})$ 
implies $\ext_R^j(S,C)_{\p}=0$ for
$j>d$.

Part~\eqref{item32} follows from~\eqref{item38}
and the definition
of a Gorenstein homomorphism.
\end{proof}

When $\vf$ is 
local,  
Theorem~\ref{lem105} shows that the next inequality can be strict
(if $\pd_R(S)>0$) or not (if $\amp(C)=0=\pd_R(S)$).
Example~\ref{ex301} shows that the converse to the final 
statement need not hold.

\begin{thm} \label{lem103}
Assume that $\vf$ is module-finite with $\fd(\vf)<\infty$, and
let $C,X$ be homologically finite $R$-complexes
with $C$ semidualizing.
There is an inequality
$$\gkdim{\rhom_R(S,C)}_S(\rhom_R(S,X))\leq
\gkdim{C}_R(X)+\amp(C).$$
In particular, if $X$ is $C$-reflexive, then $\rhom_R(S,X)$
is $\rhom_R(S,C)$-reflexive.
\end{thm}

\begin{proof}
Set $\cbc{(-)}{\vf}=\rhom_R(S,-)$. 
It suffices to verify the final statement.  Indeed, if 
$\gkdim{\cbc{C}{\vf}}_S(\cbc{X}{\vf})$ 
and 
$\gkdim{C}_R(X)$ 
are both finite, then
Theorem~\ref{lem102} 
and~\ref{para402}\eqref{item207}
explain (2) below 
\begin{align*}
\gkdim{\cbc{C}{\vf}}_S(\cbc{X}{\vf}) 
&\stackrel{(1)}{=}\inf(\cbc{C}{\vf})
-\inf(\rhom_S(\cbc{X}{\vf},\cbc{C}{\vf}))\\
&\stackrel{(2)}{\leq}
\sup(C)-\inf(S\lotimes_R\rhom_R(X,C))\\
&\stackrel{(3)}{\leq}\sup(C)-\inf(\rhom_R(X,C))\\
&\stackrel{(4)}{=}\amp(C)+\gkdim{C}_R(X)
\end{align*}
while (1) and (4) are by definition, 
and (3) follows 
from~\ref{para201}. 

Assume now that $X$ is $C$-reflexive.  The complex
$\rhom_R(X,C)$ is homologically bounded below, so~\ref{para402}\eqref{item207} 
and~\ref{para201}
imply the same for
$\rhom_S(\cbc{X}{\vf},\cbc{C}{\vf})$. 
The commutative diagram shows that the biduality morphism
$\delta^{\cbc{C}{\vf}}_{\cbc{X}{\vf}}$ is an isomorphism
$$\xymatrix{
\cbc{X}{\vf} \ar[rr]^-{\delta^{\cbc{C}{\vf}}_{\cbc{X}{\vf}}} \ar[dd]_= 
&& \rhom_S(\rhom_S(\cbc{X}{\vf},\cbc{C}{\vf}),\cbc{C}{\vf}) 
\ar[d]^{\simeq}_{\ref{para402}\eqref{item207}} \\
 && \rhom_S(S\lotimes_R\rhom_R(X,C),\rhom_R(S,C))  \\
\rhom_R(S,X) \ar[rr]^-{\rhom_R(S,\delta^C_{X})}_-{\simeq}
&& \rhom_R(S,\rhom_R(\rhom_R(X,C),C)) \ar[u]_{\simeq}^{\ref{standard}\eqref{std}} 
} 
$$ 
and it follows that $\cbc{X}{\vf}$ is $\cbc{C}{\vf}$-reflexive.
\end{proof}

When $\vf$ 
local
the first inequality in our next result
can be strict  (if $\pd_R(S)>0$)
or not (if $\pd_R(S)=0$).  We do not know 
if the second inequality 
can be strict.

\begin{thm} \label{lem105}
Let $C,X$ be homologically degreewise finite over $R$ 
with $C$ semidualizing.
If $\vf$ is module-finite with $\fd(\vf)<\infty$
and $\mspec(R)\subseteq\image(\vf^*)$, then 
\begin{align*}
\gkdim{\rhom_R(S,C)}_S(\rhom_R(S,X)) \hspace{-1cm}\\
&\leq\gkdim{C}_R(X)\\
&\leq\gkdim{\rhom_R(S,C)}_S(\rhom_R(S,X))+\pdim_R(S).
\end{align*} 
Thus, 
$\rhom_R(S,X)$
is $\rhom_R(S,C)$-reflexive if and only if $X$ is $C$-reflexive.
If either  $R$ is local or $\amp(C)=0$, 
then the second inequality 
is an equality.
\end{thm}

\begin{proof}
Set $\cbc{(-)}{\vf}=\rhom_R(S,-)$. 
First, we assume that $\cbc{X}{\vf}$ is $\cbc{C}{\vf}$-reflexive and
prove that $X$ is $C$-reflexive; the converse is in
Theorem~\ref{lem103}.
The 
complexes $\cbc{X}{\vf}$ and
$S\lotimes_R\rhom_R(X,C)$ are homologically finite
by~\ref{para402}\eqref{item207}.
Theorem~\ref{prop201}\eqref{item205} and Corollary~\ref{prop201a}\eqref{item205a}
imply the same for
$\rhom_R(X,C)$ and $X$.
In the 
commutative diagram 
from the proof of Theorem~\ref{lem103},
the morphism $\delta^{\cbc{C}{\vf}}_{\cbc{X}{\vf}}$ is an isomorphism,
hence so are $\rhom_R(S,\delta^C_{X})$  and $\delta^C_{X}$ by
Corollary~\ref{para10a}.

Now assume 
$\gkdim{C}_R(X),\gkdim{\cbc{C}{\vf}}_S(\cbc{X}{\vf})<\infty$.
The first desired inequality follows from the 
numbered sequence in the proof of 
Theorem~\ref{lem103} because $\inf(\cbc{C}{\vf})\leq\inf(C)$ by 
Theorem~\ref{lem102}. The second inequality
is in the next sequence.
\begin{align*}
\gkdim{C}_R(X)
&\stackrel{(1)}{=} \inf(C)-\inf(\rhom_R(X,C)) \\
&\stackrel{(2)}{=} 
\inf(C)-\inf(\rhom_{S}(\cbc{X}{\vf},\cbc{C}{\vf})) \\
&\stackrel{(3)}{\leq} \pdim_R(S)+\inf(\cbc{C}{\vf})
-\inf(\rhom_{S}(\cbc{X}{\vf},\cbc{C}{\vf})) \\
&\stackrel{(4)}{=} \gkdim{\cbc{C}{\vf}}_S(\cbc{X}{\vf}) +\pdim_R(S)
\end{align*}
(1) and (4) are by definition, (2) is from~\ref{para402}\eqref{item207} and 
Theorem~\ref{prop201}\eqref{item204b},
and (3) is in Theorem~\ref{lem102}. 
If $\amp(C)=0$ or $R$ is local, (3) is an equality 
by Proposition~\ref{lem104}. 
\end{proof}

Here is a version of Theorem~\ref{lem01a4} for $\rhom_R(S,-)$; its 
proof is almost identical, using Theorem~\ref{lem105} in place of 
Theorem~\ref{lem01a3}, and the isomorphism~\ref{para402}\eqref{item207}.

\begin{thm} \label{lem206}
Assume that $\vf$ is module-finite with 
$\fd(\vf)<\infty$, $\image(\vf^*)$ contains $\mspec(R)$,
and $\Pic(\vf)$ is injective.
If $C,C'$ are semidualizing $R$-complexes such that
$\rhom_R(S,C)\simeq \rhom_R(S,C')$, then $C\simeq C'$.\qed
\end{thm}

Each  inequality
in the next result
can be
strict: 
For the first and third, let $\vf$ be local and 
$\pd_R(S)>0$, and use the equality;
For the second, 
see Example~\ref{ex106}.
To see that each one can be an equality, take $\amp(C)=0=\pd_R(S)$.

\begin{prop} \label{lem201}
Let $C$ be a semidualizing $R$-complex and 
$X$ a homologically finite $S$-complex.
If $\vf$ is module-finite and $\fd(\vf)<\infty$, 
then there are inequalities
\begin{align*}
\gkdim{\rhom_R(S,C)}_S(X)-\amp(C)
&\leq\gkdim{C}_R(X)\\
&\leq\gkdim{\rhom_R(S,C)}_S(X)+\pd_R(S)
\end{align*}
with equality on the right if $R$ is local or $\amp(C)=0$.
In particular, the complex $X$ is simultaneously
$C$-reflexive and $\rhom_R(S,C)$-reflexive.
If $\image(\vf^*)$ contains $\mspec(R)$, then
$\gkdim{\rhom_R(S,C)}_S(X)
\leq\gkdim{C}_R(X)$.
\end{prop}

\begin{proof}
Set $\cbc{(-)}{\vf}=\rhom_R(S,-)$. 
Simultaneous reflexivity
is proved
in~\cite[(6.5)]{christensen:scatac}, so 
assume that $X$ is
$C$-reflexive and $\cbc{C}{\vf}$-reflexive. 
In the next sequence
\begin{align*}
\gkdim{\cbc{C}{\vf}}_S(X)
&\stackrel{(1)}{=}\inf(\cbc{C}{\vf})-\inf(\rhom_R(X,C))\\
&\stackrel{(2)}{\leq}\sup(C)-\inf(\rhom_R(X,C))\\
&\stackrel{(3)}{=}\amp(C)+\gkdim{C}_R(X)
\end{align*}
(1) is by adjunction, (2) is by Theorem~\ref{lem102}, 
and (3) is by 
definition.  This is the first inequality. 
For the second inequality, 
start with adjunction 
in (4)
\begin{align*}
\gkdim{C}_R(X)
&\stackrel{(4)}{=}\inf(C)-\inf(\rhom_S(X,\cbc{C}{\vf}))\\
&\stackrel{(5)}{\leq}\pd_R(S)+\inf(\cbc{C}{\vf})-\inf(\rhom_S(X,\cbc{C}{\vf}))\\
&\stackrel{(6)}{=}\pd_R(S)+\gkdim{\cbc{C}{\vf}}_S(X)
\end{align*}
while Theorem~\ref{lem102} 
yields (5), and (6) is by 
definition.  If $R$ is local or $\amp(C)=0$, then 
(5) is an equality by Theorem~\ref{lem102} 
and thus so is the second inequality.

If $\mspec(R)\subseteq\image(\vf^*)$, then
Corollary~\ref{prop201a} gives 
$\inf(\cbc{C}{\vf})\leq\inf(C)$. Using this in (2) above 
gives the third inequality.
\end{proof}

\begin{ex} \label{ex106}
Certain inequalities 
in Propositions~\ref{lem06b}, \ref{lem104}, and~\ref{lem201}
and in Theorems~\ref{lem01a1} and~\ref{lem102}
can be strict, even when $\spec(R)$ and $\spec(S)$ are connected.
Let $k$ be a field and set 
\begin{align*}
A&=k[X,Y,Z]/(Y^2,YZ)
&U&=A\smallsetminus((X,Y)A\cup(Y,Z)A)
&R&=U^{-1}A.
\end{align*}
The ring $R$ has two maximal ideals and one nonmaximal prime ideal 
\begin{align*}
\m&=(X,Y)R 
&\n&=(Y,Z)R
&\p&=(Y)R
\end{align*}
and $\spec(R)$ is connected as $\p\subseteq\m\cap\n$.
The minimal injective resolution of the 
normalized dualizing complex for $R$ has the form
$$D=0\to E(R/\p)\to 
E(R/\m)\oplus E(R/\n)
\to 0.$$
Since $R_{\n}$ is not Cohen-Macaulay, one has 
$\HH_1(D)_{\n}\neq 0\neq\HH_0(D)_{\n}$ 
and hence $\inf(D)=0$.
Also,  $E(R/\m)_{\p}=0=E(R/\n)_{\p}$ 
implies 
$\HH_1(D)_{\p}\cong E(R/\p)_{\p}\neq 0$. 
Since $R_{\m}$ is Gorenstein, one has $D_{\m}\sim R_{\m}$;
and since $0\neq\HH_1(D)_{\p}=(\HH_1(D)_{\m})_{\p_{\m}}$ it follows
that $\HH_1(D)_{\m}\neq 0$ and thus 
$D_{\m}\simeq\shift R_{\m}$.

With $S=R/(X)R\cong k[Y,Z]_{(Y)}/(Y^2)$ and $\vf\colon R\to S$ the natural surjection,
one has $\pd_R(S)=1$
and $\supp_R(S)=\{\m\}$.
The map $\vf_{\m S}\colon R_{\m}\to S$ is local Gorenstein
of grade 1, giving the third isomorphism below;  the other computations are routine.
\begin{gather*}
\notag
\rhom_R(S,D)
\simeq\rhom_{R_{\m}}(S,D_{\m})
\simeq\rhom_{R_{\m}}(S,\shift R_{\m})
\simeq S\\
\notag
D\lotimes_R S
\simeq D_{\m}\lotimes_{R_{\m}} S
\simeq \shift R_{\m}\lotimes_{R_{\m}} S
\simeq\shift S\\
\gkdim{D}_R(S)=0<1=
\sup\{\gkdim{D_{\m}}_{R_{\m}}(S_{\m}),\gkdim{D_{\n}}_{R_{\n}}(S_{\n})\} \\
\gkdim{D}_R(S) 
=0<1=\pd_R(S)\\
\inf(D)=0<1=\inf(D\lotimes_R S)\\
\amp(D\lotimes_R S)=0<1=\amp(D) \\
\inf(D)-\pd_R(S)=-1<0=\inf(\rhom_R(S,D))\\
\gkdim{D}_R(S)=0<1=\gkdim{\rhom_R(S,D)}_S(S)+\pd_R(S)
\end{gather*}
\end{ex}

The next result follows from Propositions~\ref{lem203} and~\ref{lem201}.
If $\vf$ is local, the inequality 
can be strict (if $\pd(S)>0$) or not (if $\pd(S)=0=\amp(C)$); see 
Theorem~\ref{lem107}. 

\begin{thm} \label{lem106}
Let $C,X$ be homologically finite $R$-complexes
with $C$ semidualizing.
When $\vf$ is module-finite with $\fd(\vf)<\infty$
there
is an inequality
$$\gkdim{\rhom_R(S,C)}_S(X\lotimes_R S)\leq
\gkdim{C}_R(X)+\amp(C)+\pd_R(S).$$
In particular, if $X$ is $C$-reflexive, then $X\lotimes_R S$
is $\rhom_R(S,C)$-reflexive.\qed
\end{thm}

Here 
is Theorem III from the introduction.
When $\vf$ is local, 
the equality guarantees that the inequalities are strict if and only if $\pd_R(S)>0$.

\begin{thm} \label{lem107}
Let $C,X$ be homologically finite $R$-complexes
with $C$ semidualizing.
If $\vf$ is module-finite with $\fd(\vf)<\infty$
and $\mspec(R)\subseteq\image(\vf^*)$, then
\begin{align*}
\gkdim{C}_R(X)-\pd_R(S)
&\leq\gkdim{\rhom_R(S,C)}_S(X\lotimes_R S)\\
&\leq\gkdim{C}_R(X)+\pd_R(S).
\end{align*}
Thus, 
$X\lotimes_R S$
is $\rhom_R(S,C)$-reflexive if and only if $X$ is $C$-reflexive.
If 
$R$ is local or
$\amp(C)=0=\amp(\rhom_R(S,R))$, then
$$\gkdim{\rhom_R(S,C)}_S(X\lotimes_R S)=\gkdim{C}_R(X).$$
\end{thm}

\begin{proof}
Set $\cbc{(-)}{\vf}=\rhom_R(S,-)$. 
In the following sequence
\begin{align*}
\gkdim{\cbc{C}{\vf}}_S(X\lotimes_R S)
&\stackrel{(1)}{\leq} \gkdim{C}_R(X\lotimes_R S)\\
&\stackrel{(2)}{\leq} \gkdim{C}_R(X)+\pd_R(S)\\
&\stackrel{(3)}{\leq} \gkdim{C}_R(X\lotimes_R S)+\pd_R(S)\\
&\stackrel{(4)}{\leq} \gkdim{\cbc{C}{\vf}}_S(X\lotimes_R S)+2\pd_R(S)
\end{align*}
(1) and (4) are in 
Proposition~\ref{lem201}, and (2) and (3) are in Theorem~\ref{lem204}.
When one of the extra conditions holds,
there is a similar sequence
$$
\gkdim{\cbc{C}{\vf}}_S(X\lotimes_R S)
=\gkdim{C}_R(X\lotimes_R S)-\pd_R(S)
=\gkdim{C}_R(X)
$$
by 
Theorem~\ref{lem204} and
Proposition~\ref{lem201}.
\end{proof}

\begin{ex} \label{ex301}
Without the hypothesis on $\mspec(R)$ in
Theorems~\ref{lem204}, 
\ref{lem01a3}, 
\ref{lem205}, 
\ref{lem105}, 
and~\ref{lem107}, the nontrivial implications fail, even when $\spec(R)$ is connected:
one can have 
$\gkdim{C}_R(X)=\infty$ even though 
each of the following is finite:
$\gkdim{C}_R(\rhom_R(S,X))$, 
$\gkdim{\rhom_R(S,C)}_R(\rhom_R(S,X))$,
$\gkdim{C}_R(X\lotimes_R S)$, 
$\gkdim{C\lotimes_R S}_R(X\lotimes_R S)$,
$\gkdim{\rhom_R(S,C)}_R(X\lotimes_R S)$.

Let 
$(R_0,\m_0)$ 
be a non-Gorenstein local ring and set $R=R_0[Y]$ with $\m=(\m_0,Y)R$ 
and
$X=R/\m\oplus R$ and $S=R/(Y-1)$ with the natural surjection
$R\to S$.  Then $X$ is not $R$-reflexive 
since if it were then $R_{\m}/\m R_{\m}$ would be $R_{\m}$-reflexive 
implying that $R_{\m}$ is Gorenstein.  However,
$X\lotimes_R S\simeq S$ has finite projective dimension over 
$S$ and over $R$, so it is reflexive with respect to each 
complex that is $R$-semidualizing or $S$-semidualizing; 
similarly for 
$\rhom_R(S,X)\simeq\shift^{-1}S$. Finally, 
$\spec(R)$ is connected as the existence of nontrivial idempotents 
in $R$ would give rise to such elements in $R_0$;  see, e.g., 
\cite[Exer.~1.22]{atiyah:ica}.
\end{ex}

\section{Factorizable 
local homomorphisms of finite flat dimension:  Cobase change}
\label{sec6}

Motivated by~\cite{avramov:rhafgd,iyengar:golh} we extend results of Section~\ref{sec7}
to special non-finite cases.

\begin{prop} \label{lem207}
Let $\dot\vf\colon R\to R'$ and $\vf'\colon R'\to S$ be 
homomorphisms of finite flat dimension
with $\vf'$ module-finite and $X$ a 
homologically degreewise
finite $R$-complex.  
\begin{enumerate}[\rm\quad(a)]
\item \label{item214}
If the $R$-complex $X$ is homologically bounded (respectively, 
semidualizing), then 
the $S$-complex $\rhom_{R'}(S,X\lotimes_R R')$ 
is so as well. 
\item \label{item215}
Assume that $\dot\vf$ is faithfully flat and $\image((\vf')^*)$
contains 
$\mspec(R')$.   
If the $S$-complex $\rhom_{R'}(S,X\lotimes_R R')$
is homologically bounded (respectively, semidualizing), then 
the $R$-complex $X$ is so as well.
\end{enumerate}
\end{prop}

\begin{proof}
\eqref{item214} If $X$ is homologically bounded, then so is 
$\rhom_{R'}(S,X\lotimes_R R')$ by~\ref{para201} 
and~\ref{para201a}.  
Theorems~\ref{lem01a1} 
and~\ref{lem102} 
yield the other implication.

\eqref{item215} When 
$\dot{\vf}$ is flat, the  
isomorphism
$\HH_i(X\lotimes_R R')\cong \HH_i(X)\otimes_R R'$
implies that the $R'$-complex
$X\lotimes_R R'$ is homologically degreewise finite.  
If $\rhom_{R'}(S,X\lotimes_R R')$ is homologically bounded, then
so is $X\lotimes_R R'$ by Corollary~\ref{prop201a}\eqref{item205a},
and so is $X$.
Theorems~\ref{lem01a1} 
and~\ref{lem102} 
provide the remaining implication.
\end{proof}

For the rest of this paper, we focus on local homomorphisms that factor nicely.

\begin{para} \label{factor}
When $\vf$ is local, a \emph{regular} 
(respectively, \emph{Gorenstein}) 
\emph{factorization} of $\vf$ is a pair of local homomorphisms
$R\xra{\dot\vf}R'\xra{\vf'}S$
such that $\vf=\vf'\dot\vf$, $\vf'$ is surjective, and $\dot{\vf}$ is flat with 
regular (respectively, Gorenstein) closed fibre.  
In either case, the homomorphisms $\vf$ and $\vf'$ have finite flat 
dimension simultaneously by~\cite[(3.2)]{foxby:daafuc}.
When the ring
$R'$ 
is complete, the regular factorzation is a \emph{Cohen 
factorization}.  
It is straightforward to construct a regular factorization 
when $\vf$ is essentially of finite type.
Also, if $S$ is complete, then
$\vf$ 
admits a Cohen factorization~\cite[(1.1)]{avramov:solh}. 
\end{para}

\begin{lem} \label{lemGF}
Assume that $\vf$ is module-finite and local and that it admits a
Gorenstein factorization $R\xra{\dot\vf_1}R_1\xra{\vf_1'}S$ 
with $R,R_1,S$  complete.
Then there exists a commutative diagram of local homomorphisms
\[
\xymatrix{
 & R_1\ar@{>>}[d]_{\pi}\ar@{>>}[rd]^{\vf_1'} \\
R \ar[ur]^{\dot\vf_1}\ar[r]^{\ddot{\vf}} & R'' \ar@{>>}[r]^{\vf''} & S
} \]
where $\pi$ is surjective with kernel generated by an 
$R_1$-sequence 
and the bottom row is a Gorenstein 
factorization of $\vf$ such that $\ddot{\vf}$ is module-finite.  
\end{lem}

\begin{proof}
Since $\vf$ is module finite, 
the closed fibre $S/\m S\cong R_1/(\Ker(\vf_1'),\m)$ is Artinian and
the extension of
residue fields $k\to l$ is finite.
In particular, the ideal $(\Ker(\vf_1'),\m)R_1/\m R_1$ is primary to the maximal 
ideal of $R_1/\m R_1$.  Let $\y=y_1,\ldots,y_d\in\Ker(\vf_1')$ be a system 
of parameters for $R_1/\m R_1$, that is, a maximal $R_1/\m R_1$-sequence.  
Set $R''=R_1/(\y)$ with natural surjection $\pi\colon R_1\to R''$, and 
let the maps $\ddot{\vf}\colon R\to R''$ and 
$\vf''\colon R''\to S$
be induced by $\dot\vf_1$ and $\vf_1'$, respectively.

One has $\vf''\ddot{\vf}=\vf_1'\dot\vf_1=\vf$,
and $\vf''$ is surjective because $\vf_1'$ is 
so.  The closed fibre of $\ddot{\vf}$ is 
$R''/\m R''\cong (R_1/\m R_1)/(\y)$ which is Gorenstein because 
$R_1/\m R_1$ is so.
The sequence $\y$ is $R_1$-regular,
and the map $\ddot{\vf}$ is flat;  see, e.g., \cite[Corollary to (22.5)]{matsumura:crt}.  
Finally, the equality in the next sequence is straightforward
\[ \len_R(R''/\m R'')=\len_{R''/\m R''}(R''/\m 
R'')\cdot\rank_{k}(l)<\infty.\]
and the inequality is by construction.
So, 
$\vf''$ is module-finite by~\cite[(8.4)]{matsumura:crt}.  
\end{proof}

\begin{prop} \label{prop601}
Assume 
that $\vf$ is local
and admits
Gorenstein factorizations 
$R\xra{\dot{\vf}_1}R_1\xra{\vf_1'}S$ and
$R\xra{\dot{\vf}_2}R_2\xra{\vf_2'}S$ with each $R_i$ complete.
There exists
a commutative diagram of local ring homomorphisms
\[ \xymatrix{
 & R_1 \ar@{>>}[rd]^{\vf_1'} \\
R \ar[ur]^{\dot{\vf}_1}\ar[r]^{\dot{\vf}}\ar[rd]_{\dot{\vf}_2}
& R' \ar@{>>}[u]_-{\pi_1}\ar@{>>}[d]^-{\pi_2}\ar@{>>}[r]^-{\vf'} & S \\
 & R_2 \ar@{>>}[ur]_{\vf_2'}
} \]
where $\vf'\dot{\vf}$ is a Cohen factorization of $\vf$ and each 
$\pi_i$ is surjective and Gorenstein.
\end{prop}

\begin{proof}
Taking Cohen factorizations $R\xra{\ddot\vf_i}R_i'\xra{\vf_i''}R_i$ of $\dot\vf_i$,
it is evident that the diagrams $R\xra{\ddot\vf_i}R_i'\xra{\vf_i'\vf_i''}S$
are Cohen factorizations of $\vf$.  Since $\dot\vf_i$ is flat with Gorenstein closed
fibre, the surjection $\vf_i''$ is Gorenstein by~\cite[(2.4)]{avramov:lgh}
and~\cite[(3.2)]{avramov:solh}. 
The `comparison theorem' for Cohen 
factorizations~\cite[(1.2)]{avramov:solh}, 
provides a commutative diagram of local ring homomorphisms
\[ \xymatrix{
 & R_1' \ar@{>>}[rd]^{\vf_1'\vf_1''} \\
R \ar[ur]^{\ddot{\vf}_1}\ar[r]^{\dot{\vf}}\ar[rd]_{\ddot{\vf}_2}
& R' \ar@{>>}[u]_-{\tau_1}\ar@{>>}[d]^-{\tau_2}\ar@{>>}[r]^{\vf'} & S \\
 & R_2' \ar@{>>}[ur]_{\vf_2'\vf_2''}
} \]
where $\vf'\dot{\vf}$ is a Cohen factorization of $\vf$ and each 
$\tau_i$ is surjective with kernel generated by a regular sequence.  
Each $\tau_i$ is Gorenstein by~\cite[(4.3)]{avramov:lgh}, and hence
so is each $\pi_i=\vf_i''\tau_i$.  Thus, these maps yield a 
diagram with the desired properties.
\end{proof}

\begin{thm} \label{thm5}
Let $X$ be a homologically finite $R$-complex.
Assume that $\vf$ is local with $\fd(\vf)$ finite and that $\vf$ admits
Gorenstein factorizations 
$R\xra{\dot{\vf}_1}R_1\xra{\vf_1'}S$ and
$R\xra{\dot{\vf}_2}R_2\xra{\vf_2'}S$.
Set 
$d=\depth(\vf)$ and $d_i=\depth(\dot{\vf}_i)$ for $i=1,2$.
\begin{enumerate}[\quad\rm(a)]
\item \label{itemA}
The $S$-complexes $\shift^{d_1}\rhom_{R_1}(S,X\lotimes_R R_1)$ and 
$\shift^{d_2}\rhom_{R_2}(S,X\lotimes_R R_2)$ are isomorphic. 
\item \label{itemC}
When $\vf$ is Gorenstein at $\n$, the $S$-complexes
$\shift^{d_i}\rhom_{R_i}(S,X\lotimes_R R_i)$ and
$\shift^{d}X\lotimes_R S$
are isomorphic.
\item \label{itemB}
When $\vf$ is module-finite, the $S$-complexes
$\shift^{d_i}\rhom_{R_i}(S,X\lotimes_R R_i)$ and
$\rhom_R(S,X)$ are isomorphic.
\end{enumerate}
\end{thm}

\begin{proof}
First, we show that, if $\vf$ is module-finite and Gorenstein at $\n$,
then the $S$-complexes $\shift^{d}X\lotimes_R S$ and $\rhom_R(S,X)$
are isomorphic.  To this end, 
note that $\grade_R(S)=-d$ and so Proposition~\ref{prop01a}\eqref{item32}
implies
$\shift^{d}S\simeq\rhom_R(S,R)$
since $S$ is local.
This provides the first of 
the following isomorphisms
$$\shift^{d}S\lotimes_R X
\simeq \rhom_R(S,R)\lotimes_R X
\simeq\rhom_R(S,X)$$
where the other is from~\ref{para202}\eqref{item217}.  
This establishes the desired 
isomorphism.

The completed diagrams 
$\comp{R}\xra{\comp{\dot\vf_i}}\comp{R_i}\xra{\comp{\vf_i'}}\comp{S}$
are Gorenstein factorizations of $\comp{\vf}\colon \comp{R}\to\comp{S}$.
Using Lemma~\ref{sri}, one can replace the given factorizations with the completed ones to assume that
the local rings $R,R_1,R_2, S$ are complete.

By considering the upper and lower halves of the
diagram provided by Proposition~\ref{prop601} 
we assume that there is a commutative diagram
of local homomorphisms 
$$\xymatrix{
 & R_1 \ar@{>>}[d]_{\tau}\ar@{>>}[rd]^{\vf_1'} \\
R \ar[ur]^{\dot\vf_1}\ar[r]^{\dot\vf_2} & R_2 \ar@{>>}[r]^{\vf_2'} & S
} 
$$ 
where $\tau$ is surjective and Gorenstein. 
By definition then, one has $d_2=d_1+\depth(\tau)$.

\eqref{itemA}  The above diagram gives a sequence of 
isomorphisms
\begin{align*}
\shift^{d_2}\rhom_{R_2}(S,X\lotimes_R R_2) \hspace{-2cm} \\
& \stackrel{(1)}{\simeq} 
\shift^{d_2}\rhom_{R_2}(S,(X\lotimes_R R_1)\lotimes_{R_1}R_2) \\
& \stackrel{(2)}{\simeq} 
\shift^{d_2}\rhom_{R_2}(S,\shift^{-\depth(\tau)}\rhom_{R_1}(R_2,X\lotimes_R R_1)) \\
& \stackrel{(3)}{\simeq} 
\shift^{d_1}\rhom_{R_2}(S,\rhom_{R_1}(R_2,X\lotimes_R R_1)) \\
& \stackrel{(4)}{\simeq} 
\shift^{d_1}\rhom_{R_1}(S,X\lotimes_R R_1) 
\end{align*}
where (1) is by associativity, (2) follows from the
the first paragraph since 
$\tau$ is Gorenstein and surjective, (3) follows from the 
final observation of the previous paragraph, and (4) is adjunction. 

\eqref{itemC}  When $\vf$ is Gorenstein, the same is true of each 
$\vf_i'$ by~\cite[(3.2)]{avramov:solh} and~\cite[(2.4)]{avramov:lgh}.  
Since each $\vf_i'$ is also surjective, the first paragraph 
gives the first isomorphism in the next sequence
where the second isomorphism is associativity and cancellation.
\[ \shift^{d_i}\rhom_{R_i}(S,X\lotimes_R R_i)
\simeq \shift^{d_i+\depth(\vf_i')}(X\lotimes_R R_i)\lotimes_{R_i}S 
\simeq \shift^{d}X\lotimes_R S \]

\eqref{itemB}  When $\vf$ is module-finite, 
the diagram provided by Lemma~\ref{lemGF} 
yields a sequence of isomorphisms where $d''=\depth(\ddot\vf)$.
\begin{align*}
\shift^{d_1}\rhom_{R_1}(S,X\lotimes_R R_1) 
& \stackrel{(5)}{\simeq} 
\shift^{d''}\rhom_{R''}(S,X\lotimes_R R'') \\
& \stackrel{(6)}{\simeq} 
\shift^{d''}\rhom_{R''}(S, 
\shift^{-d''} \rhom_R(R'',X)) \\
& \stackrel{(7)}{\simeq} 
\rhom_R(S,X)
\end{align*}
(5) is by part~\eqref{itemA},
(6)
follow from the first paragraph, 
and (7) is by adjunction.
\end{proof}

We employ the following handy notation for the remainder of this 
section.

\begin{para} \label{cbc}
Assume that $\vf$ is local with $\fd(\vf)$ finite and admits
a Gorenstein factorization $R\xra{\dot{\vf}}R'\xra{\vf'}S$
with $d=\depth(\dot{\vf})$.
For a homologically finite $R$-complex $X$, set  
\[ \cbc{X}{\vf}=\shift^{d}\rhom_{R'}(S,X\lotimes_R R'). \]
Theorem~\ref{thm5} shows that this is 
independent of the choice of Gorenstein factorization and that
$\cbc{X}{\vf}\simeq\rhom_R(S,X)$ when $\vf$ is module-finite.
\end{para}

\begin{disc} \label{cbcD}
With $\vf$ as in~\ref{cbc}, the complex $\cbc{R}{\vf}$ is normalized dualizing for $\vf$.
If $D$ is a (normalized)
dualizing complex for $R$, then the complex 
$\cbc{D}{\vf}$ is (normalized) dualizing for $S$;  see 
Proposition~\ref{prop7}. 
\end{disc}

Next is an alternate description of $\cbc{X}{\vf}$
that follows directly from~\ref{para202}\eqref{item217}.
In it, we tensor over $S$ in order to stress that
complexes are isomorphic over $S$
and not just over $R$.  
A similar remark applies to Proposition~\ref{prop502}.

\begin{prop} \label{prop501}
If $\vf$ is as in~\ref{cbc} and $X$ is a homologically finite $R$-complex,
then there is an isomorphism
$\cbc{X}{\vf}\simeq(X\lotimes_R S)\lotimes_S\cbc{R}{\vf}$. 
\qed
\end{prop}

The next isomorphisms follows from parts~\eqref{item207} and~\eqref{item208} 
of~\ref{para402}.

\begin{prop} \label{prop502}
If $\vf$ is as in~\ref{cbc} 
then there are isomorphisms
\begin{gather*}
\rhom_S(\cbc{X}{\vf},\cbc{Y}{\vf})
\simeq \rhom_R(X,Y)\lotimes_R S \\
\rhom_S(X\lotimes_R S,\cbc{Y}{\vf})
\simeq (\rhom_R(X,Y)\lotimes_R S)\lotimes_S \cbc{R}{\vf} 
\end{gather*}
for all homologically finite $R$-complexes  $X,Y$. \qed
\end{prop}

\begin{prop} \label{prop7}
Assume that $\vf$ is local with $\fd(\vf)$ finite 
and let $C$ be a semidualizing
$R$-complex.
The Poincar\'{e} and Bass series of $C\lotimes_R S$ are
\begin{align*}
P^S_{C\lotimes_R S}(t)&=P^R_C(t) &
I_S^{C\lotimes_R S}(t)&=I_R^C(t)I_{\vf}(t). \end{align*}
If $\vf$ has a Gorenstein factorization, then the 
Poincar\'{e} and Bass series of $\cbc{C}{\vf}$ are
\begin{align*}
P^S_{\cbc{C}{\vf}}(t)&=P^R_C(t)I_{\vf}(t) & 
I_S^{\cbc{C}{\vf}}(t)&=I_R^C(t).\end{align*}
\end{prop}

\begin{proof}
The first Poincar\'{e} series
is from~\cite[(1.5.3)]{avramov:rhafgd}, and the Bass series follows  
\[ I_S^{C\lotimes_R S}(t)
\stackrel{(1)}{=}I^S_S(t)/P^S_{C\lotimes_R S}(t) 
\stackrel{(2)}{=}I^R_R(t)I_{\vf}(t)/P^R_C(t)
\stackrel{(3)}{=}I_R^C(t)I_{\vf}(t) \]
where (1) and (3) are by~\ref{para101} and (2) is
from~\ref{paraBass}.
If $\vf$ admits a Gorenstein factorization, then the second Poincar\'{e}
series follows from Proposition~\ref{prop501} with~\cite[(1.5.3)]{avramov:rhafgd}
and~\cite[(1.7.6)]{christensen:scatac}, and
the second Bass series is computed like the first one.
\end{proof}

Here we record the analogue of Theorem~\ref{lem105} for our new setting.

\begin{cor} \label{thm2}
If $\vf$ is as in~\ref{cbc} and $C,X$ are homologically finite 
$R$-complexes with $C$ semidualizing, then
$\gkdim{\cbc{C}{\vf}}_{S}(\cbc{X}{\vf})=\gkdim{C}_R(X) +\depth(\vf)$. 
In particular, $\cbc{X}{\vf}$ is $\cbc{C}{\vf}$-reflexive if and only if 
$X$ is 
$C$-reflexive.
\end{cor}

\begin{proof}
Let $R\xra{\dot\vf}R'\xra{\vf'}S$ be a Gorenstein factorization of $\vf$ 
and set $d=\depth(\dot\vf)$.
Equalities (1) and (5) in the following sequence are by 
definition
\begin{align*}
\gkdim{\cbc{C}{\vf}}_{S}(\cbc{X}{\vf})
&\stackrel{(1)}{=}\gkdim{\shift^d\cbc{(C\lotimes_R 
R')}{\vf'}}_{S}(\shift^d\cbc{(X\lotimes_R R')}{\vf'}) \\
&\stackrel{(2)}{=}\gkdim{\cbc{(C\lotimes_R 
R')}{\vf'}}_{S}(\cbc{(X\lotimes_R R')}{\vf'})+d \\
&\stackrel{(3)}{=}\gkdim{C\lotimes_R 
R'}_{R'}(X\lotimes_R R')-\pd_{R'}(S)+d \\
&\stackrel{(4)}{=}\gkdim{C}_{R}(X)+\depth(\vf')+d \\
&\stackrel{(5)}{=}\gkdim{C}_{R}(X)+\depth(\vf)
\end{align*}
while (2) is by~\cite[(3.12)]{christensen:scatac},
(3) is Theorem~\ref{lem105}, and
(4) is from Theorem~\ref{lem01a3} and the Auslander-Buchsbaum formula.
\end{proof}

Theorems~\ref{lem01a4} and~\ref{lem206} provide the proof of the next 
result.

\begin{cor} \label{lem208}
Let $\vf$ be as in~\ref{cbc}.  When $C,C'$ are
semidualizing $R$-complexes, one has  
$\cbc{C}{\vf}\simeq  \cbc{C'}{\vf}$ if and 
only if $C\simeq C'$.\qed
\end{cor}

Replace Proposition~\ref{lem201} with Theorem~\ref{lem107} in the proof 
of Corollary~\ref{thm2} to prove the next result.

\begin{cor} \label{thm4}
If $\vf$ is as in~\ref{cbc} and $C,X$ are homologically finite 
$R$-complexes with $C$ semidualizing, then
$\gkdim{\cbc{C}{\vf}}_{S}(X\lotimes_R S)=\gkdim{C}_R(X)$. 
In particular, $X\lotimes_R S$ is $\cbc{C}{\vf}$-reflexive 
if and only if $X$ is 
$C$-reflexive.\qed
\end{cor}

\begin{disc}
With the  reflexivity relations of Theorem~\ref{lem01a3}
and Corollaries~\ref{thm2} and~\ref{thm4} in mind,
we 
wish 
to characterize
the finiteness of 
$\gkdim{C\lotimes_R S}_S(\cbc{X}{\vf})$.  If $\gkdim{C}_R(X)$ is 
finite and $\vf$ is Gorenstein at $\n$, then $\gkdim{C\lotimes_R 
S}_S(\cbc{X}{\vf})$ is finite by 
Theorems~\ref{lem01a3}
and~\ref{thm5}\eqref{itemC}. 
We wonder if the converse holds.  Here is one 
instance of this:
If $\gkdim{C\lotimes_R S}_S(\cbc{C}{\vf})$ is finite, then  
$C\lotimes_R S$ and $\cbc{C}{\vf}$ are shift isomorphic 
by Lemma~\ref{prop202},
and~\cite[(3.7(c))]{frankild:sdcms} 
implies that $\vf$ is 
Gorenstein at $\n$.
\end{disc}

\begin{prop} \label{propCBC}
Let $\vf$ be local with $\fd(\vf)$ finite 
and $C,C'$ semidualizing
$R$-complexes such that $C'$ is $C$-reflexive.
There are coefficientwise
equalities
\begin{align*}
P^S_{\rhom_S(C'\lotimes_R S,C\lotimes_R S)}(t) 
& =P^R_{\rhom_R(C',C)}(t)   \\
I_S^{\rhom_S(C'\lotimes_R S,C\lotimes_R S)}(t) 
&=I_R^{\rhom_R(C',C)}(t)I_{\vf}(t). 
\end{align*}
If $\vf$ has a Gorenstein factorization, then there are equalities
\begin{align*}
P^S_{\rhom_S(\cbc{C'}{\vf},\cbc{C}{\vf})}(t)
& =P^R_{\rhom_R(C',C)}(t) \\
I_S^{\rhom_S(\cbc{C'}{\vf},\cbc{C}{\vf})}(t)
& =I_R^{\rhom_R(C',C)}(t)I_{\vf}(t)   \\
P^S_{\rhom_S(C'\lotimes_R S,\cbc{C}{\vf})}(t) 
& =P^R_{\rhom_R(C',C)}(t)I_{\vf}(t)  \\
I_S^{\rhom_S(C'\lotimes_R S,\cbc{C}{\vf})}(t) 
&=I_R^{\rhom_R(C',C)}(t).
\end{align*}
\end{prop}

\begin{proof}
For the first Poincar\'{e} series, use~\cite[(1.7.6)]{christensen:scatac}
with~\ref{para402}\eqref{item206}.
When $\vf$ admits 
a Gorenstein factorization, the other Poincar\'{e} series 
come from Proposition~\ref{prop502} with~\cite[(1.5.3)]{avramov:rhafgd}
and~\cite[(1.7.6)]{christensen:scatac}.
The Bass series follow 
as in
Proposition~\ref{prop7}.
\end{proof}

\section*{Acknowledgments}
 
A.F.~is grateful to the Department of Mathematics at the 
University of Illinois at Urbana-Champaign for its hospitality while 
much of this research was conducted.  S.S.-W.~is similarly 
grateful to the Institute for Mathematical 
Sciences at the University of Copenhagen.
Both authors express their gratitude to L.~Avramov and
L.~W.~Christensen
for 
stimulating conversations and helpful comments about this 
research, to S.~Iyengar for allowing us to include
Lemma~\ref{sri}, and to the anonymous referee for improving
the presentation.

\providecommand{\bysame}{\leavevmode\hbox to3em{\hrulefill}\thinspace}
\providecommand{\MR}{\relax\ifhmode\unskip\space\fi MR }
\providecommand{\MRhref}[2]{
  \href{http://www.ams.org/mathscinet-getitem?mr=#1}{#2}
}
\providecommand{\href}[2]{#2}

\end{document}